\documentclass[onecolumn, 12pt, english]{IEEEtran}

\usepackage[T1]{fontenc}
\usepackage[latin9]{inputenc}
\usepackage{color}
\usepackage{array}
\usepackage{float}
\usepackage{amsmath}
\usepackage{amsthm}
\usepackage{amssymb}
\usepackage{graphicx}

\makeatletter

\floatstyle{ruled}
\newfloat{algorithm}{tbp}{loa}
\providecommand{\algorithmname}{Algorithm}
\floatname{algorithm}{\protect\algorithmname}

\theoremstyle{plain}
\newtheorem{thm}{\protect\theoremname}
\theoremstyle{plain}
\newtheorem{lem}[thm]{\protect\lemmaname}

\theoremstyle{example}

\theoremstyle{definition}

\theoremstyle{plain}

\theoremstyle{plain}

\usepackage{babel}
\providecommand{\lemmaname}{Lemma}
\providecommand{\theoremname}{Theorem}

\makeatother

\usepackage{babel}
\providecommand{\lemmaname}{Lemma}
\providecommand{\theoremname}{Theorem}

\begin{document}

\title{Distributed Zeroth-Order Stochastic Optimization in Time-varying Networks}

\author{Wenjie~Li,~\IEEEmembership{Member, IEEE}, and Mohamad~Assaad,~\IEEEmembership{Senior Member, IEEE}
\thanks{W. Li was with  Laboratoire des Signaux et Systèmes,
CentraleSupélec, 91190 Gif-sur-Yvette, France. He is now with
Paris Research Center, Huawei Technologies, 92100 Boulogne-Billancourt, France (e-mail: liwenjie28@huawei.com) }\thanks{M. Assaad is with Laboratoire des Signaux et
Systèmes (L2S, UMR CNRS 8506), CentraleSupélec, 91190 Gif-sur-Yvette, France (e-mail: mohamad.assaad@centralesupelec.fr).}}
\maketitle
\begin{abstract}
  We consider a distributed convex optimization problem in a network
which is time-varying and not always strongly connected. The local
cost function of each node is affected by some stochastic process.
All nodes of the network collaborate to minimize the average of their
local cost functions. The major challenge of our work is that the
gradient of cost functions is supposed to be unavailable and has to
be estimated only based on the numerical observation of cost functions.
Such problem is known as zeroth-order stochastic convex optimization
(ZOSCO). In this paper we take a first step towards the distributed
optimization problem with a ZOSCO setting. The proposed algorithm
contains two basic steps at each iteration: i) each unit updates a
local variable according to a random perturbation based single point
gradient estimator of its own local cost function; ii) each unit exchange
its local variable with its direct neighbors and then perform a weighted
average. In the situation where the cost function is smooth and strongly
convex, our attainable optimization error is $O(T^{-1/2})$ after
$T$ iterations. This result is interesting as $O(T^{-1/2})$ is the
optimal convergence rate in the ZOSCO problem. We have also investigate
the optimization error with the general Lipschitz convex function,
the result is $O(T^{-1/4})$.
\end{abstract}

% REQUIRED

\section{Introduction }

Distributed optimization is an attractive problem in various domain.
A set of distributed nodes need to minimize some global cost function,
while each node can only have access to information of local cost
function. For example, in the Internet of Things applications, each
device collects local data to serve the whole network. Most of the
existing work in this area requires the knowledge of gradient information.

In this work, we consider the distributed optimization problem with
more realistic yet challenging settings: 1) the network is time-varying
and not always strongly connected; 2) the local cost functions are
disturbed by some stochastic process; 3) the gradient information
is unavailable and each node is only able to access to a numerical
value of its local function at each time. Such problem belongs to
zeroth-order stochastic convex optimization (ZOSCO), it can be also
named as derivative-free stochastic convex optimization or bandit
convex optimization. Using the numerical value of objective function,
one can only has some \textit{biased} estimation of its gradient. It is thus
not surprising that there is a fundamental performance gap between
ZOSCO and Stochastic Approximation (SA) \cite{nemirovski2009robust}, as in SA  an \textit{unbiased} gradient estimation is available. Several
lower bounds of the convergence rate have been derived in \cite{jamieson2012query,shamir2013complexity,duchi2015optimal}
to state that, without the knowledge of gradient, the optimization
error cannot be better than $O\left(T^{-0.5}\right)$ after $T$ iterations.
Note that when the gradient is available, the convergence rate can
be as fast as $O\left(T^{-1}\right)$ given that the objective function
is strongly concave \cite{nemirovski2009robust}. 
Recently, ZOSCO is receiving an increasing interest in various applications. For example,
\cite{chen2018bandit} considered the management of fog computing
in IoT, i.e., nodes in the fog layer need to collaboratively process
some data requests to minimize the processing delay. Due to the unpredictable
network congestion, the closed form expression of delay can not be
available. As a result, one can only apply ZOSCO. An example related
to sensor selection for parameter estimation was considered in \cite{liu2017zeroth}.
The main goal is to find the optimal tradeoff between sensor activations
and the estimation accuracy. In order to evaluate the gradient of
the objective function, the calculation of inverse matrices is necessary, 
while a zeroth-order method proposed \cite{liu2017zeroth} helps to
avoid matrix inversion. Zeroth-order methods are also interesting
for solving adversarial machine learning problem, see \cite{liu2018zeroth2}
as an example.

Although both distributed optimization and ZOSCO are popular problems
in recent years, the distributed optimization problem in a ZOSCO setting
is rarely investigated to the best our knowledge. In this work, we
make a first step towards such problem. Similar to the typical distributed
optimization algorithm, our algorithm contains a gradient step and
a communication step in each iteration. While the major difference
is that, in our gradient step, each node has to estimate the gradient
of its local function by using a \textit{single} numerical value of the function. 

We consider a time-varying and not always strongly connected network
and derive the convergence rate of the proposed algorithm. We have
obtained a $O\left(T^{-0.5}\right)$ convergence rate of the optimization
error, under the assumption that the global function is strongly convex
and each local function is smooth. This result is remarkable as
it is the same as the lower bound of the convergence rate in the centralized
ZOSCO problem. We have also considered the case with general Lipschitz
convex function (non smooth and non strongly convex), the resulted
attainable convergence rate is $O\left(T^{-0.25}\right)$. While similar
assumptions related to network graph have been considered in \cite{nedic2009distributed},
 the major difference is that gradient was supposed to be available
in \cite{nedic2009distributed}. The main challenge in our work
comes from the fact that our zeroth-order approach leads to \textit{biased}
estimation of the gradient and therefore we have to analyze an additional
term associated to the gradient estimation bias. 

\textbf{Related work.} There exist significant amount of work related
to distributed optimization aiming at improving convergence rate,
while most of the work consider the direct availability of gradient
information, e.g., \cite{nedic2009distributed,seaman2017optimal,scaman2018optimal}.
ZOSCO or bandit convex optimization problem has also attracted much
interests recently. Several works have based on two-point gradient estimator
(TPGE) \cite{agarwal2010optimal,duchi2015optimal,liu2018zeroth,hajinezhad2019zone,nesterov2017},
which requires two successive observation of the cost function $F(\boldsymbol{\theta}_{t};\boldsymbol{\xi}_{t})$
and $F(\boldsymbol{\theta}_{t}';\boldsymbol{\xi}_{t})$ under the
same stochastic parameter $\boldsymbol{\xi}_{t}$. Several recent
work \cite{yuan2015zeroth,sahu2018distributed,pang2019randomized,hajinezhad2019zone}
have applied such TPGE to address distributed zeroth-order convex
optimization. However, TPGE cannot be realized in the situation where
the value of $\boldsymbol{\xi}_{t}$ change fast (e.g., in an i.i.d.
manner), as it is impossible to observe both $F(\boldsymbol{\theta}_{t};\boldsymbol{\xi}_{t})$
and $F(\boldsymbol{\theta}_{t}';\boldsymbol{\xi}_{t})$. Such a situation arises widely in various applications in practice, e.g. in wireless networks. In this work, we focus on single point gradient estimator, which has
been barely addressed in distributed ZOSCO problem. In fact, it is possible to have a biased estimation of
gradient by using a single realization of the cost function \cite{flaxman2005online}.
In \cite{shamir2013complexity}, the gradient estimator is shown to
be unbiased for quadratic functions, while it is a biased estimation for general nonlinear function. In the framework of ZOSCO, there exist several work that
provide advanced methods trying to achieve the optimal convergence
rate for general convex functions \cite{hazan2014bandit,hazan2016optimal,bubeck2017kernel},
however, within a centralized setting as the proposed algorithms contains operations of vectors and matrices which requires to be handled by a central controller.

In our previous work \cite{li2021IT,li2018cdc, li2021TAC}, we considered a different framework
where each node controls its own action to perform
ZOSCO. The network model was also different: each node has a constant
probability to communicate with any other node of the network. However,
in this work, each node aim to optimize a common vector and we consider
a time-varying communication matrix between the nodes.

The rest of the paper is organized as follows. Section~\ref{sec:Problem-setting-and}
describes the problem setting and the basic assumptions. Section~\ref{sec:algo}
presents the proposed distributed optimization algorithm. In Section~\ref{sec:smooth_strong},
we provide the convergence rate of the proposed algorithm for smooth
and strongly-convex function. Whereas general Lipschitz convex function
is considered in Section~\ref{sec:general}. Section~\ref{sec:Conclusion}
concludes this paper.

\section{Problem setting and assumptions \label{sec:Problem-setting-and}}

In this paper, matrices are in boldface upper case letters and vectors
are in boldface lower case letters. Calligraphic font denotes set.
We denote by $\left\Vert \boldsymbol{x}\right\Vert $ the Euclidean
norm of any vector$\boldsymbol{x}$. $\mathbf{I}_{M}$ denotes a identity
matrix of size $M\times M$ and $\mathbf{0}_{M}$ denotes a zero matrix
of size $M\times M$.

\subsection{Problem setting}

Consider a set of nodes $\mathcal{V}=\left\{ 1,2,\ldots,N\right\} $
and a time-varying directed graph $\mathcal{G}_{t}=\left(\mathcal{V},\mathcal{E}_{t}\right)$
at any timeslot $t$. The edge set $\mathcal{E}_{t}$ denotes a collection
of pairs having direct communication. We can introduce a time-varying
communication matrix $\mathbf{A}\left(t\right)=\left[A_{i,j}\left(t\right)\right]_{i,j\in\mathcal{V}}$,
with $A_{i,j}\left(t\right)\in\left(0,1\right)$ only if $\left(i,j\right)\in\mathcal{E}_{t}$,
otherwise $A_{i,j}\left(t\right)=0$, $\forall\left(i,j\right)\notin\mathcal{E}_{t}$. 

For any $\boldsymbol{\theta}\in\mathcal{K}\subseteq\mathbb{R}^{M}$
and some random vector $\boldsymbol{\xi}_{i}\in\mathcal{S}_{i}$,
denote $F_{i}\left(\boldsymbol{\theta};\boldsymbol{\xi}_{i}\right):\mathcal{K}\times\mathcal{S}_{i}\rightarrow\mathbb{R}$
as the local cost function of node~$i$. Note that we use $\boldsymbol{\xi}_{i}$
to describe some non-additive stochastic process, which is assumed
to be i.i.d. and ergodic. Introduce the expected local cost function
$f_{i}\left(\boldsymbol{\theta}\right)=\mathbb{E}_{\boldsymbol{\xi}}\left[F_{i}\left(\boldsymbol{\theta};\boldsymbol{\xi}_{i}\right)\right]$.
The problem is to find the optimal value of $\boldsymbol{\theta}$
to minimize the average cost function, i.e.,
\begin{equation}
\min_{\boldsymbol{\theta}\in\mathcal{K}}f\left(\boldsymbol{\theta}\right)=\frac{1}{N}\sum_{i=1}^{N}f_{i}\left(\boldsymbol{\theta}\right)=\frac{1}{N}\sum_{i=1}^{N}\mathbb{E}_{\boldsymbol{\xi}}\left[F_{i}\left(\boldsymbol{\theta};\boldsymbol{\xi}_{i}\right)\right],\label{eq:prob_1}
\end{equation}
under the challenging condition that:
\begin{enumerate}
\item each node~$i$ only knows the numerical value of $F_{i}$ rather
then the closed form expression of $F_{i}$;
\item each node~$i$ can only communicate with its neighbors at each time,
the global information of the network is not available.
\end{enumerate}
We denote the minimizer of $f\left(\boldsymbol{\theta}\right)$ as
$\boldsymbol{\theta}^{*}$, i.e., 
\begin{equation}
\boldsymbol{\theta}^{*}=\arg\min_{\boldsymbol{\theta}\in\mathcal{K}}f\left(\boldsymbol{\theta}\right)=\frac{1}{N}\sum_{i=1}^{N}f_{i}\left(\boldsymbol{\theta}\right).
\end{equation}

A typical strategy to solve the distributed optimization problem (\ref{eq:prob_1})
is to make each node\ $i$ update a local variable $\boldsymbol{\theta}_{i}$
aiming to minimize its local cost function, then nodes need communication
to have an agreement of their local variables. An alternative way
to write problem (\ref{eq:prob_1}) is the following

\begin{align}
\min_{\boldsymbol{\theta}_{1},\ldots,\boldsymbol{\theta}_{N}\in\mathcal{K}} & \frac{1}{N}\sum_{i=1}^{N}f_{i}\left(\boldsymbol{\theta}\right)=\frac{1}{N}\sum_{i=1}^{N}\mathbb{E}_{\boldsymbol{\xi}}\left[F_{i}\left(\boldsymbol{\theta}_{i};\boldsymbol{\xi}_{i}\right)\right],\label{eq:problem}\\
\mathrm{s.t.} & \boldsymbol{\theta}_{i}=\boldsymbol{\theta}_{2}=\ldots=\boldsymbol{\theta}_{N}.\label{eq:equal}
\end{align}

\subsection{Assumptions}

In this section we present the assumptions that will be considered
in this paper. In terms of the communication matrix, we assume that:

\textbf{A1}: i). $\mathbf{A}\left(t\right)$ is doubly-stochastic,
i.e., $\sum_{j\in\mathcal{N}}A_{ij}\left(t\right)=\sum_{j\in\mathcal{N}}A_{ji}\left(t\right)=1$,
$\forall i\in\mathcal{N}$; ii). every non-zero element of $\mathbf{A}\left(t\right)$
is always equal or larger than a constant $a>0$; iii). two arbitrary
nodes can directly communicate with each other at least once every
$\tau>0$ consecutive time slots, in other words, $(\mathcal{V},\bigcup_{t=t_{0}}^{t_{0}+\tau-1}\mathcal{E}_{t})$
is strongly connected for any $t_{0}\geq0$.

Throughout this paper, we have the following assumptions in terms
of the basic properties of cost function.

\textbf{A2}: The feasible set $\mathcal{K}$ is bounded, convex, and
contains the zero point $\mathbf{0}$. For any $\boldsymbol{\theta}\in\mathcal{K}$,
we have $\left\Vert \boldsymbol{\theta}\right\Vert \leq R<+\infty$. 

\textbf{A3}: The value of $F_{i}$ is bounded, we have $\left|F_{i}\left(\boldsymbol{\theta};\boldsymbol{\xi}_{i}\right)\right|\leq C<+\infty$
for any $i\in\mathcal{V}$, $\boldsymbol{\theta}\in\mathcal{K}$,
and $\boldsymbol{\xi}_{i}\in\mathcal{S}_{i}$.

Specifically, in Section~\ref{sec:smooth_strong}, we perform the
analysis assuming that the cost functions are smooth and strongly
convex as described in what follows.

\textbf{A4:} The expected global cost function $f\left(\boldsymbol{\theta}\right)$
is $\mu-$strongly convex. For any $i\in\mathcal{V}$, $f_{i}\left(\boldsymbol{\theta}\right)$
is $L_{i}$-smooth, i.e.,
\[
\left\Vert \nabla f_{i}\left(\boldsymbol{\theta}\right)-\nabla f_{i}\left(\boldsymbol{\theta}'\right)\right\Vert \leq L_{i}\left\Vert \boldsymbol{\theta}-\boldsymbol{\theta}'\right\Vert ,\qquad\forall\boldsymbol{\theta},\boldsymbol{\theta}'\in\mathcal{K},\forall i\in\mathcal{V},
\]
denote $L=\frac{1}{N}\sum_{i=1}^{N}L_{i}$.

We have also considered a more general situation in Section~\ref{sec:general}
that each local cost function is Lipschitz:

\textbf{A5:} For any $i\in\mathcal{V}$, $f_{i}\left(\boldsymbol{\theta}\right)$
is $\ell_{i}$-smooth, i.e.,
\[
\left\Vert f_{i}\left(\boldsymbol{\theta}\right)-f_{i}\left(\boldsymbol{\theta}'\right)\right\Vert \leq\ell_{i}\left\Vert \boldsymbol{\theta}-\boldsymbol{\theta}'\right\Vert ,\qquad\forall\boldsymbol{\theta},\boldsymbol{\theta}'\in\mathcal{K},\forall i\in\mathcal{V}.
\]

\section{Distributed zeroth-order stochastic optimization algorithm\label{sec:algo}}

In this section, we first introduce a general single point gradient
estimator and then present a simple algorithm to solve the distributed
optimization problem as well as its fundamental properties. 

\subsection{A general single point gradient estimator}

We consider a gradient estimator which is based on a single point
evaluation of local cost function. At each time~$t$, node~$i$
estimates the gradient of its local cost function with
\begin{equation}
\widehat{\boldsymbol{g}}_{i}\left(t\right)=\frac{\boldsymbol{\nu}_{i,t}F_{i}\left(\boldsymbol{\theta}_{i}\left(t\right)+\beta_{t}\boldsymbol{\nu}_{i,t};\boldsymbol{\xi}_{i,t}\right)}{\beta_{t}},\label{eq:g_est}
\end{equation}
where $\beta_{t}>0$ represents the step-size and $\boldsymbol{\nu}_{i,t}\in\left[-V,V\right]^{M}\subseteq\mathbb{R}^{M}$
is a random perturbation vector which is independently generated by
each node. Two different settings will be considered in Sections~\ref{sec:smooth_strong}
and~\ref{sec:general} respectively. 

\textbf{S1}: For any $i\in\mathcal{V}$ and $t\in\left\{ 1,\ldots,T\right\} $,
$\left\Vert \boldsymbol{\nu}_{i,t}\right\Vert \leq\sqrt{M}V<+\infty$,
each element of $\boldsymbol{\nu}_{i,t}$ is i.i.d. such that $\mathbb{E}[\boldsymbol{\nu}_{i,t}]=\mathbf{0}$
and $\mathbb{E}[\boldsymbol{\nu}_{i,t}\boldsymbol{\nu}_{i,t}^{\mathsf{Tr}}]=\mathbf{I}_{M}$.

\textbf{S2}: For any $i\in\mathcal{V}$ and $t\in\left\{ 1,\ldots,T\right\} $,
$\boldsymbol{\nu}_{i,t}$ is a random unit vector such that $\mathbb{E}[\boldsymbol{\nu}_{i,t}]=\mathbf{0}$
and $\left\Vert \boldsymbol{\nu}_{i,t}\right\Vert =1$.

Clearly, S1 is a relaxed setting compared with S2: $\boldsymbol{\nu}_{i,t}$
has to be a unit vector in S2, which is not the case in S1. Providing
a single point gradient estimator under S1 is an interesting extension
of ZOSCO, which can be useful also in some other frameworks where
the perturbation vector has to be generated in a distributed manner,
i.e., each distributed node only generates some coordinate of the
vector and the resulted entire perturbation vector cannot be unit
without a control center.

In what follows, we can show that $\widehat{\boldsymbol{g}}_{i}\left(t\right)$
is a reasonable estimator of $\nabla f_{i}$ with either S1 or S2,
as stated in Lemma~\ref{lem:bias_bound} and Lemma~\ref{lem:bias_general}.
\begin{lem}
\label{lem:bias_bound}Suppose that the random perturbation vector
$\boldsymbol{\nu}_{i,t}$ is generated according to S1, then under
Assumption A4, we can derive that 
\begin{equation}
\nabla f_{i}\left(\boldsymbol{\theta}_{i}\left(t\right)\right)=\mathbb{E}\left[\widehat{\boldsymbol{g}}_{i}\left(t\right)\mid\boldsymbol{\theta}_{i}\left(t\right)\right]-\boldsymbol{b}_{i}\left(t\right),
\end{equation}
where $\boldsymbol{b}_{i}\left(t\right)$ represents a estimation
bias of the gradient estimator with the following property 
\begin{equation}
\left\Vert \boldsymbol{b}_{i}\left(t\right)\right\Vert \leq\frac{1}{2}M^{\frac{3}{2}}V^{2}L_{i}\beta_{t}=O\left(\beta_{t}\right)\qquad\forall i\in\mathcal{V},t\in\mathbb{N}^{+}\label{eq:lemma_1}
\end{equation}
\end{lem}
Lemma~\ref{lem:bias_bound} can be proved by applying Taylor's theorem
and mean value theorem, the proof detail can be found in Appendix~\ref{subsec:appx_a}.
In general, the gradient estimator $\widehat{\boldsymbol{g}}_{i}\left(t\right)$
introduces an non-zero bias term $\boldsymbol{b}_{i}\left(t\right)$.
Nevertheless, we can see from Lemma~\ref{lem:bias_bound} that $\left\Vert \boldsymbol{b}_{i}\left(t\right)\right\Vert $
decreases with the step-size $\beta_{t}$, if $\beta_{t}$ is vanishing.
Lemma~\ref{lem:bias_bound} will be a basis of our analysis concerning
smooth and strongly convex functions in Section~\ref{sec:smooth_strong}. 
\begin{lem}
\label{lem:bias_general}Suppose that $\boldsymbol{\nu}_{i,t}$ is
generated according to S2, introduce a smoothed version of function
$f_{i}$ for any $i\in\mathcal{N}$, i.e., 
\[
\widetilde{f}_{i}\left(\boldsymbol{\theta}\right)=\mathbb{E}_{\boldsymbol{\varpi}\in\mathbb{R}^{M}:\left\Vert \boldsymbol{\varpi}\right\Vert \leq1}\left[f_{i}\left(\boldsymbol{\theta}+\beta\boldsymbol{\varpi}\right)\right],
\]
then we have 
\[
\frac{\beta}{M}\nabla\widetilde{f}_{i}\left(\boldsymbol{\theta}\right)=\mathbb{E}_{\boldsymbol{\nu}_{i,t}\in\in\mathbb{R}^{M}:\left\Vert \boldsymbol{\nu}_{i,t}\right\Vert =1}\left[f_{i}\left(\boldsymbol{\theta}+\beta\boldsymbol{\nu}_{i,t}\right)\boldsymbol{\nu}_{i,t}\right]=\mathbb{E}\left[\widehat{\boldsymbol{g}}_{i}\left(t\right)\mid\boldsymbol{\theta}\right].
\]
\end{lem}
Lemma~2 is a fundamental result presented in \cite{flaxman2005online},
it states that the gradient estimator $\widehat{\boldsymbol{g}}_{i}\left(t\right)$
is in fact the exact gradient of a smoothed version of $f_{i}$ when
$\boldsymbol{\nu}_{i,t}$ satisfies S2. Our analysis related to general
Lipschitz convex functions will rely on Lemma~\ref{lem:bias_general}.

\subsection{Distributed Algorithm}

By using the gradient estimator $\widehat{\boldsymbol{g}}_{i}\left(t\right)$
as defined in (\ref{eq:g_est}), at time~$t$ each node~$i$ needs
to observe the value of $F_{i}$ at the point $\boldsymbol{\theta}_{i}\left(t\right)+\beta_{t}\boldsymbol{\nu}_{i,t}$.
In order to have $\boldsymbol{\theta}_{i}\left(t\right)+\beta_{t}\boldsymbol{\nu}_{i,t}\in\mathcal{K}$,
the value of $\boldsymbol{\theta}_{i}\left(t\right)$ should belong
to the set 
\[
\mathcal{K}_{t}=\left\{ \boldsymbol{x}\in\mathbb{R}^{M}:\boldsymbol{x}+\beta_{t}V\mathbf{1}\in\mathcal{K}\right\} ,
\]
where $\mathbf{1}$ denotes a vector of ones with dimension $M\times1$.
Recall that $\boldsymbol{\nu}_{i,t}\in\left[-V,V\right]^{M}$. A simple
solution is to use the projection operator defined as 
\begin{equation}
\Pi_{\mathcal{K}_{t}}\left(\boldsymbol{x}\right)=\arg\min_{\boldsymbol{x}'\in\mathcal{K}_{t}}\left\Vert \boldsymbol{x}-\boldsymbol{x}'\right\Vert ,\qquad\forall\boldsymbol{x}\in\mathbb{R}^{M}.\label{eq:projection}
\end{equation}
at each time~$t$. The distributed algorithm considered in this paper
is given by
\begin{align}
\boldsymbol{\theta}_{i}\left(t+1\right) & =\Pi_{\mathcal{K}_{t+1}}\left(\sum_{j=1}^{N}A_{i,j}\left(t\right)\boldsymbol{\theta}_{j}\left(t\right)-\alpha_{t}\widehat{\boldsymbol{g}}_{i}\left(t\right)\right),\label{eq:algo}
\end{align}
where $\alpha_{t}>0$ is another step-size. In order to make the algorithm
converge with a fast speed, both $\alpha_{t}$ and $\beta_{t}$ should
be chosen wisely. We will provide a detailed discussion in the next
sections. The basic idea of (\ref{eq:algo}) is to make each node
update $\boldsymbol{\theta}_{i}\left(t\right)$ according to the gradient
estimator of its local cost function, as well as the weighted sum
of $\boldsymbol{\theta}_{j}\left(t\right)$ from neighbors. The detailed
algorithm is described in Algorithm~\ref{alg:D-zosco}.

\begin{algorithm}[h]
\caption{\label{alg:D-zosco}D-ZOSCO Algorithm for each node~$i$}
\begin{enumerate}
\item Initialize $t=1$ and $\boldsymbol{\theta}_{i}\left(t\right)=\mathbf{0}$. 
\item Generate a random perturbation vector $\boldsymbol{\nu}_{i,t}$ and
observe the value of local cost function $F_{i}$ at point $\boldsymbol{\theta}_{i}\left(t\right)+\beta_{t}\boldsymbol{\nu}_{i,t}$.
\item Update the gradient estimator $\widehat{\boldsymbol{g}}_{i}\left(t\right)=\boldsymbol{\nu}_{i,t}F_{i}(\boldsymbol{\theta}_{i}\left(t\right)+\beta_{t}\boldsymbol{\nu}_{i,t};\boldsymbol{\xi}_{i,t})/\beta_{t}$.
\item Communicate with all neighbors, i.e., broadcast the value of $\boldsymbol{\theta}_{i}\left(t\right)$
and receive $\boldsymbol{\theta}_{j}\left(t\right)$ from each neighboring
node~$j$ (with $A_{i,j}\left(t\right)>0$). Compute the weighted
average $\sum_{j=1}^{N}A_{i,j}\left(t\right)\boldsymbol{\theta}_{j}\left(t\right)$.
\item Update $\boldsymbol{\theta}_{i}\left(t+1\right)$ according to (\ref{eq:algo}),
i.e., $\boldsymbol{\theta}_{i}\left(t+1\right)=\Pi_{\mathcal{K}_{t+1}}(\sum_{j=1}^{N}A_{i,j}\left(t\right)\boldsymbol{\theta}_{j}\left(t\right)-\alpha_{t}\widehat{\boldsymbol{g}}_{i}\left(t\right))$.
\item $t=t+1$, go to 2.
\end{enumerate}
\end{algorithm}

We introduce an auxiliary variable
\[
\overline{\boldsymbol{\theta}}\left(t\right)=\frac{1}{N}\sum_{i=1}^{N}\boldsymbol{\theta}_{i}\left(t\right),
\]
which is the average of local variables $\boldsymbol{\theta}_{i}\left(t\right)$
at time~$t$. Since an essential objective of our distributed algorithm
is to make $\boldsymbol{\theta}_{1}\left(t\right)\approx\ldots\approx\boldsymbol{\theta}_{N}\left(t\right)$
according to the problem setting (\ref{eq:equal}), it is desirable
to have $\boldsymbol{\theta}_{i}\left(t\right)\rightarrow\overline{\boldsymbol{\theta}}\left(t\right)$
for any $i\in\mathcal{N}$ as $t$ is large.. The following lemma
presents an upper bound of $\left\Vert \boldsymbol{\theta}_{i}\left(t\right)-\overline{\boldsymbol{\theta}}\left(t\right)\right\Vert $,
from which we can see the disagreement of local variables $\boldsymbol{\theta}_{i}\left(t\right)$
from a quantitative point of view.
\begin{lem}
\label{lem:dist_aver}Suppose that Assumptions~A1, A2, A3 hold, besides
$\left\Vert \boldsymbol{\nu}_{i,t}\right\Vert \leq\sqrt{M}V$ and
the step-sizes are vanishing such that $\alpha_{t}=\alpha_{0}t^{-c_{1}}$,
$\beta_{t}=\beta_{0}t^{-c_{2}}$, and $\frac{\alpha_{t}}{\beta_{t}}=\frac{\alpha_{0}}{\beta_{0}}t^{-c_{3}}$,
then we have $\left\Vert \boldsymbol{\theta}_{i}\left(t\right)-\overline{\boldsymbol{\theta}}\left(t\right)\right\Vert \leq\delta_{t}$
with 
\begin{equation}
\delta_{t}=\left(2^{1+c_{3}}+N\rho\left(\frac{1+c_{3}}{1-\eta}+\frac{c_{3}}{\left(1-\eta\right)^{2}}\right)\right)\sqrt{M}VC\frac{\alpha_{t}}{\beta_{t}},\label{eq:delta_def}
\end{equation}
where $\rho=2\frac{1+a^{\left(N-1\right)\tau}}{1-a^{\left(N-1\right)\tau}}$
and $\eta=\left(1-a^{\left(N-1\right)\tau}\right)^{\frac{1}{\left(N-1\right)\tau}}$.
If $\left\Vert \boldsymbol{\nu}_{i,t}\right\Vert =1$ and the step-sizes
are constant such that $\alpha_{t}=\alpha$ and $\beta_{t}=\beta$,
then we have
\begin{equation}
\left\Vert \boldsymbol{\theta}_{i}\left(t\right)-\overline{\boldsymbol{\theta}}\left(t\right)\right\Vert \leq\delta=\left(2+N\rho\frac{\eta}{1-\eta}\right)C\frac{\alpha}{\beta}.\label{eq:delta_def-1}
\end{equation}
\end{lem}
The proof detail is provided in Appendix~\ref{subsec:appx_b}. 

\section{Convergence rate for smooth strongly convex function\label{sec:smooth_strong}}

\subsection{Main result}

In this section, we present the convergence rate of average optimization
error for a special class of convex function. 

Denote $\overline{\boldsymbol{\theta}}^{\left(T\right)}=\sum_{t=1}^{T}\sum_{i=1}^{N}\boldsymbol{\theta}_{i}\left(t\right)/T$.
The aim of this section is to find an upper bound for the average
optimization error, defined as $\mathbb{E}[f(\overline{\boldsymbol{\theta}}^{\left(T\right)})-f(\boldsymbol{\theta}^{*})]$. 
\begin{thm}
\label{thm:result_strong}Suppose that Assumptions~A1-A4 hold, $\boldsymbol{\theta}^{*}\in\mathcal{K}_{0}$
and each node~$i$ updates $\boldsymbol{\theta}_{i}\left(t\right)$
according to (\ref{eq:algo}) with perturbation vector $\boldsymbol{\nu}_{i,t}$
satisfying S1 and step-sizes $\alpha_{t}=\frac{3}{\mu}t^{-1}$ and
$\beta_{t}=\left(6\lambda_{\mathsf{I}}^{-2}\lambda_{\mathsf{III}}\right)^{\frac{1}{4}}t^{-\frac{1}{4}}$,
then 
\begin{align}
 & \mathbb{E}\left[f\left(\overline{\boldsymbol{\theta}}^{\left(T\right)}\right)-f\left(\boldsymbol{\theta}^{*}\right)\right]\nonumber \\
 & \leq\frac{L}{N}\Psi^{*}\lambda_{\mathsf{I}}^{2}\left(\frac{1}{\sqrt{T}}+\frac{\lambda_{\mathsf{II}}}{\mu\sqrt{6\lambda_{\mathsf{III}}}}\frac{\log\left(T\right)}{T}+\left(\frac{\lambda_{\mathsf{II}}}{\mu\sqrt{6\lambda_{\mathsf{III}}}}+\frac{\lambda_{\mathsf{II}}^{2}}{4\mu^{2}\lambda_{\mathsf{III}}}\right)\frac{1}{T}\right)\label{eq:rate_strong-1}\\
 & =O\left(M^{2}T^{-\frac{1}{2}}\right),\nonumber 
\end{align}
where the constant terms are defined as
\[
\begin{cases}
\lambda_{\mathsf{I}}=\sqrt{NM^{3}}V^{2}L\\
\lambda_{\mathsf{II}}=2L\sqrt{NM}VC\left(2^{\frac{7}{4}}+N\rho\frac{10-7\eta}{4\left(1-\eta\right)^{2}}\right)\\
\lambda_{\mathsf{III}}=2NMV^{2}C^{2}\left(\left(2^{\frac{7}{4}}+N\rho\frac{10-7\eta}{4\left(1-\eta\right)^{2}}\right)^{2}+1\right)\\
\Psi^{*}=\max\left\{ 4\sqrt{6}NR^{2}\left(\lambda_{\mathsf{I}}+\frac{\lambda_{\mathsf{I}}\lambda_{\mathsf{II}}}{6\mu\sqrt{\lambda_{\mathsf{III}}}}\right)^{-2},\frac{\sqrt{6\lambda_{\mathsf{III}}}}{\lambda_{\mathsf{I}}\mu^{2}}\right\} 
\end{cases}
\]
\end{thm}
Note that we introduce the constant terms $\lambda_{\mathsf{I}}$,
$\lambda_{\mathsf{II}}$, and $\lambda_{\mathsf{III}}$ mainly to
lighten the notations. We explain the basic steps of our analysis
in the rest of this section.

\subsection{Proof sketch of Theorem~4}

Thanks to the convexity and smoothness of the cost functions, we can
show that (with proof detail in Appendix~\ref{subsec:appx_c})
\begin{equation}
\mathbb{E}\left[f\left(\overline{\boldsymbol{\theta}}^{\left(T\right)}\right)-f\left(\boldsymbol{\theta}^{*}\right)\right]\leq\frac{L}{2NT}\sum_{t=1}^{T}\mathbb{E}\left[\sum_{i=1}^{N}\left\Vert \boldsymbol{\theta}_{i}\left(t\right)-\boldsymbol{\theta}^{*}\right\Vert ^{2}\right]\label{eq:error_divergence}
\end{equation}
from which we find that the convergence speed of $\mathbb{E}[f(\overline{\boldsymbol{\theta}}^{\left(T\right)})-f(\boldsymbol{\theta}^{*})]$
is the same as that of the average divergence defined as 
\begin{equation}
d_{t}=\mathbb{E}\left[\sum_{i=1}^{N}\left\Vert \boldsymbol{\theta}_{i}\left(t\right)-\boldsymbol{\theta}^{*}\right\Vert ^{2}\right].\label{eq:aver_dist_def}
\end{equation}
Thus we can focus on the upper bound of $d_{t}$. Our first step is
to find the relation between $d_{t+1}$ and $d_{t}$, as described
in the following lemma.
\begin{lem}
\label{lem:algo_iter_strong}Suppose that Assumptions~A1-A4 hold,
$\boldsymbol{\theta}^{*}\in\mathcal{K}_{0}$, and each node~$i$
updates $\boldsymbol{\theta}_{i}\left(t\right)$ according to (\ref{eq:algo}),
and $\alpha_{t}=\alpha_{0}t^{-1}$, $\beta_{t}=\beta_{0}t^{-\frac{1}{4}}$.
Consider the constants $\lambda_{\mathsf{I}}$, $\lambda_{\mathsf{II}}$,
and $\lambda_{\mathsf{III}}$ as defined in Theorem~\ref{thm:result_strong},
we have
\begin{equation}
d_{t+1}\leq\left(1-2\mu\alpha_{t}\right)d_{t}+\left(\lambda_{\mathsf{I}}+\lambda_{\mathsf{II}}\frac{\alpha_{t}}{\beta_{t}^{2}}\right)\alpha_{t}\beta_{t}\sqrt{d_{t}}+\lambda_{\mathsf{III}}\frac{\alpha_{t}^{2}}{\beta_{t}^{2}},\label{eq:evo_strong}
\end{equation}
\end{lem}
The proof of Lemma~\ref{lem:algo_iter_strong} uses the results presented
in Lemma~\ref{lem:bias_bound} and Lemma~\ref{lem:dist_aver}, we
need to investigate both $\left\Vert \overline{\boldsymbol{\theta}}\left(t\right)-\boldsymbol{\theta}_{i}\left(t\right)\right\Vert $
and $\left\Vert \overline{\boldsymbol{\theta}}\left(t\right)-\boldsymbol{\theta}^{*}\right\Vert $.
Please refer to Appendix~\ref{subsec:appx_d} for the proof details. 

Based on Lemma~5, the next step is to show the upper bound of $d_{t}$
by using induction.
\begin{lem}
\label{lem:rate_strong}Consider $\alpha_{t}=\alpha_{0}t^{-1}$ and
$\beta_{t}=\beta_{0}t^{-\frac{1}{4}}$. Suppose that Assumptions A1-A6
hold, $\boldsymbol{\theta}^{*}\in\mathcal{K}_{0}$, and $\alpha_{0}>\frac{3}{4\mu}$.
Introduce $t_{0}=\left\lceil 2\mu\alpha_{0}\right\rceil $, then
\begin{equation}
d_{t}\leq\Psi\left(\lambda_{\mathsf{I}}+\lambda_{\mathsf{II}}\frac{\alpha_{0}}{\beta_{0}^{2}}t^{-\frac{1}{2}}\right)^{2}t^{-\frac{1}{2}},\qquad\forall t\geq1\label{eq:rate_strong}
\end{equation}
with the constant term
\begin{equation}
\Psi=\max\left\{ \frac{4NR^{2}t_{0}^{\frac{1}{2}}}{\left(\lambda_{\mathsf{I}}+\lambda_{\mathsf{II}}\frac{\alpha_{0}}{\beta_{0}^{2}}t_{0}^{-\frac{1}{2}}\right)^{2}},\left(\frac{\alpha_{0}\beta_{0}+\sqrt{\alpha_{0}^{2}\beta_{0}^{2}+2\left(4\mu\alpha_{0}-3\right)\frac{\lambda_{\mathsf{III}}\alpha_{0}^{2}}{\lambda_{\mathsf{I}}^{2}\beta_{0}^{2}}}}{4\mu\alpha_{0}-3}\right)^{2}\right\} \label{eq:const_strong}
\end{equation}
\end{lem}
The proof of Lemma~6 can be found in Appendix~\ref{subsec:appx_e}.
We find that it is possible to have $d_{t}=O\left(t^{-\frac{1}{2}}\right)$.
The constant term $\Psi$ is a function of $\alpha_{0}$ and $\beta_{0}$
with complicated form, it is necessary to properly choose $\alpha_{0}$
and $\beta_{0}$ to make $\Psi$ as tight as possible. The following
lemma provides a reasonable solution. The proof detail is provided
in Appendix~\ref{subsec:appx_f}. 
\begin{lem}
\label{lem:step_const}If $\alpha_{0}>\frac{3}{4\mu}$ and $\beta_{0}>0$,
we have
\[
\frac{\alpha_{0}\beta_{0}+\sqrt{\alpha_{0}^{2}\beta_{0}^{2}+2\left(4\mu\alpha_{0}-3\right)\frac{\lambda_{\mathsf{III}}\alpha_{0}^{2}}{\lambda_{\mathsf{I}}^{2}\beta_{0}^{2}}}}{4\mu\alpha_{0}-3}\geq\left(6\lambda_{\mathsf{I}}^{-2}\lambda_{\mathsf{III}}\right)^{\frac{1}{4}}\mu^{-1}
\]
with equality as $\alpha_{0}^{*}=\frac{3}{\mu}$ and $\beta_{0}^{*}=\left(6\lambda_{\mathsf{I}}^{-2}\lambda_{\mathsf{III}}\right)^{\frac{1}{4}}$.
\end{lem}
We can introduce the values of $\alpha_{0}^{*}$ and $\beta_{0}^{*}$
into Lemma~\ref{lem:rate_strong} to get 
\[
d_{t}\leq\Psi^{*}\left(\lambda_{\mathsf{I}}+\frac{3\lambda_{\mathsf{I}}\lambda_{\mathsf{II}}}{\mu\sqrt{6\lambda_{\mathsf{III}}}}t^{-\frac{1}{2}}\right)^{2}t^{-\frac{1}{2}}
\]
with $\Psi^{*}$ is defined in Theorem~\ref{thm:result_strong}. 

Finally, we can turn to evaluate the average optimization error $\mathbb{E}[f(\overline{\boldsymbol{\theta}}^{\left(T\right)})-f(\boldsymbol{\theta}^{*})]$
by using (\ref{eq:error_divergence}) and the upper bound of $d_{t}$,
i.e.,

\begin{align}
 & \mathbb{E}\left[f\left(\overline{\boldsymbol{\theta}}^{\left(T\right)}\right)-f\left(\boldsymbol{\theta}^{*}\right)\right]\nonumber \\
 & \leq\frac{L}{2NT}\sum_{t=1}^{T}d_{t}\leq\frac{L}{2NT}\Psi^{*}\lambda_{\mathsf{I}}^{2}\sum_{t=1}^{T}\left(t^{-\frac{1}{2}}+\frac{\sqrt{6}\lambda_{\mathsf{II}}}{\mu\sqrt{\lambda_{\mathsf{III}}}}t^{-1}+\frac{3\lambda_{\mathsf{II}}^{2}}{2\mu^{2}\lambda_{\mathsf{III}}}t^{-\frac{3}{2}}\right)\nonumber \\
 & \leq\frac{L}{N}\Psi^{*}\lambda_{\mathsf{I}}^{2}\left(\frac{1}{\sqrt{T}}+\frac{\sqrt{6}\lambda_{\mathsf{II}}}{2\mu\sqrt{\lambda_{\mathsf{III}}}}\frac{\log\left(T\right)}{T}+\left(\frac{\sqrt{6}\lambda_{\mathsf{II}}}{2\mu\sqrt{\lambda_{\mathsf{III}}}}+\frac{9\lambda_{\mathsf{II}}^{2}}{4\mu^{2}\lambda_{\mathsf{III}}}\right)\frac{1}{T}\right),\label{eq:strong_rate}
\end{align}
which is obtained by using the bounds $\sum_{t=1}^{T}t^{-\frac{1}{2}}<\int_{x=0}^{T}x^{-\frac{1}{2}}dx=2\sqrt{T}$,
$\sum_{t=1}^{T}t^{-1}<1+\int_{x=1}^{T}x^{-1}dx=\log\left(T\right)+1$,
and $\sum_{t=1}^{T}t^{-\frac{3}{2}}<1+\int_{x=1}^{\infty}x^{-\frac{3}{2}}dx=3$.

\subsection{Discussion}

In this section, we have shown that the average optimization error
after $T$ iterations can be $O(T^{-1/2})$ for smooth and strongly
convex function. Such result is nice as the optimal convergence rate
for a centralized ZOSCO problem is also $O(T^{-1/2}$). Recall that
the challenge of our problem not only comes from the ZOSCO setting,
but is also caused by the imperfect network topology: in a distributed
setting, a node only has the local information of the global function.
As derived in Lemma~\ref{lem:bias_bound} and Lemma~\ref{lem:dist_aver},
the estimation bias of gradient decreases with $\beta_{t}=O(t^{-1/4})$,
while $\left\Vert \overline{\boldsymbol{\theta}}\left(t\right)-\boldsymbol{\theta}^{*}\right\Vert $
caused by communication decreases with a faster speed $\frac{\alpha_{t}}{\beta_{t}}=O(t^{-3/4})$.
This explains that the convergence rate of our problem is more sensitive
to the performance of the biased gradient estimator.

\section{Convergence rate for Lipschitz convex function \label{sec:general}}

In this section, we consider a more general case where $f$ is neither
smooth nor strongly convex. In the situation where the cost functions
are non-smooth, we have to use the random unit perturbation vector
to ensure that the gradient estimator works well. 

Our main result is stated as follows.
\begin{thm}
\label{thm:result_convex}Suppose that Assumptions~A1, A2, A3, A5
hold, and the perturbation vector satisfy S2. If each node~$i$ updates
$\boldsymbol{\theta}_{i}\left(t\right)$ according to (\ref{eq:algo})
with constant step-sizes\textup{
\begin{equation}
\begin{cases}
\alpha_{t}=\alpha^{*}=\sqrt{\frac{2\sqrt{2}MR^{3}}{C\ell\sqrt{7+5N\rho\frac{\eta}{1-\eta}+\left(N\rho\frac{\eta}{1-\eta}\right)^{2}}}}T^{-\frac{3}{4}}\\
\beta_{t}=\beta^{*}=\sqrt{\frac{MRC}{\ell}\sqrt{14+10N\rho\frac{\eta}{1-\eta}+2\left(N\rho\frac{\eta}{1-\eta}\right)^{2}}}T^{-\frac{1}{4}}
\end{cases}\label{eq:ab_opt}
\end{equation}
}and $\boldsymbol{\theta}^{*}\in\widetilde{\mathcal{K}}=\left\{ \boldsymbol{\theta}\in\mathcal{K}:\boldsymbol{\theta}+\beta\boldsymbol{v}\in\mathcal{K},\forall\boldsymbol{v}\textrm{ s.t. }\left\Vert \boldsymbol{v}\right\Vert =1\right\} $,
then the average optimization error is bounded by \textup{
\begin{align}
 & \mathbb{E}\left[f\left(\overline{\boldsymbol{\theta}}^{\left(T\right)}\right)-f\left(\boldsymbol{\theta}^{*}\right)\right]\nonumber \\
 & \leq4\sqrt{MRC\ell\sqrt{14+10N\rho\frac{\eta}{1-\eta}+2\left(N\rho\frac{\eta}{1-\eta}\right)^{2}}}T^{-\frac{1}{4}}+\sqrt{2}R\ell T^{-\frac{1}{2}}\label{eq:general_result}\\
 & =O\left(\sqrt{M}T^{-\frac{1}{4}}\right)+O\left(T^{-\frac{1}{2}}\right)
\end{align}
}
\end{thm}
The proof sketch of Theorem~\ref{thm:result_convex} is presented
in the rest of this section.

\subsection{Proof sketch of Theorem~8}

The proof of Theorem~\ref{thm:result_convex} is based on the application
of Lemma~\ref{lem:bias_general}, which states that our gradient
estimator is the exact gradient of a smoothed version of the objective
function 
\[
\widetilde{f}\left(\boldsymbol{\theta}\right)=\frac{1}{N}\sum_{i\in\mathcal{N}}\widetilde{f}_{i}\left(\boldsymbol{\theta}\right)=\frac{1}{N}\sum_{i\in\mathcal{N}}\mathbb{E}_{\boldsymbol{\varpi}\in\mathbb{R}^{M}:\left\Vert \boldsymbol{\varpi}\right\Vert \leq1}\left[f_{i}\left(\boldsymbol{\theta}+\beta\boldsymbol{\varpi}\right)\right]
\]
 with $\boldsymbol{\theta}\in\widetilde{\mathcal{K}}\subseteq\mathcal{K}$.
It is worth mentioning that $\widetilde{\mathcal{K}}$ is a subset
of $\mathcal{K}$ to ensure $\boldsymbol{\theta}+\beta\boldsymbol{\varpi}\in\mathcal{K}$,
i.e., $\widetilde{\mathcal{K}}=\left\{ \boldsymbol{\theta}\in\mathcal{K}:\boldsymbol{\theta}+\beta\boldsymbol{v}\in\mathcal{K},\forall\boldsymbol{v}\textrm{ s.t. }\left\Vert \boldsymbol{v}\right\Vert =1\right\} $.
Introduce $\widetilde{\boldsymbol{\theta}}^{*}$ as the minimizer
of $\widetilde{f}\left(\boldsymbol{\theta}\right)$, i.e.,
\[
\widetilde{\boldsymbol{\theta}}^{*}=\arg\min_{\boldsymbol{\theta}\in\widetilde{\mathcal{K}}}\frac{1}{N}\sum_{i\in\mathcal{N}}\widetilde{f}_{i}\left(\boldsymbol{\theta}\right).
\]

First, we can build a relation between the optimization error $\mathbb{E}[f(\overline{\boldsymbol{\theta}}^{\left(T\right)})-f(\boldsymbol{\theta}^{*})]$
of the original cost function and that of the smoothed function, as
stated in the following lemma.
\begin{lem}
\label{lem:general_lem1}Suppose that Assumptions~A1, A2, A3, and
A5 hold. Each node~$i$ updates $\boldsymbol{\theta}_{i}\left(t\right)$
according to (\ref{eq:algo}) with perturbation vector satisfying
S2 and constant step-sizes $\alpha_{t}=\alpha$ and $\beta_{t}=\beta$.
Then we have
\begin{align}
\mathbb{E}\left[f\left(\overline{\boldsymbol{\theta}}^{\left(T\right)}\right)-f\left(\boldsymbol{\theta}^{*}\right)\right] & \leq\frac{1}{NT}\sum_{t=1}^{T}\sum_{i=1}^{N}\left(\mathbb{E}\left[\widetilde{f}_{i}\left(\boldsymbol{\theta}_{i}\left(t\right)\right)\right]-\widetilde{f}_{i}\left(\widetilde{\boldsymbol{\theta}}^{*}\right)\right)+2\ell\beta+\ell\delta.\label{eq:genral_error-1}
\end{align}
recall that $\delta\propto\alpha/\beta$ is defined in (\ref{eq:delta_def-1}).
\end{lem}
The gap of the two optimization errors mainly comes from the differences
between $\widetilde{f}\left(\boldsymbol{\theta}\right)$ and $f\left(\boldsymbol{\theta}\right)$,
as well as the difference between $\widetilde{f}(\widetilde{\boldsymbol{\theta}}^{*})$
and $f\left(\boldsymbol{\theta}^{*}\right)$. Please refer to the
proof details of Lemma~\ref{lem:general_lem1} in Appendix~\ref{subsec:appx_g}. 

Clearly, our next step is to find the upper bound of the optimization
error of the smoothed function $\widetilde{f}$.
\begin{lem}
\label{lem:general_lem2}Under the same condition of Lemma~\ref{lem:general_lem1},
we have
\begin{align}
\frac{1}{NT}\sum_{t=1}^{T}\sum_{i=1}^{N}\left(\mathbb{E}\left[\widetilde{f}_{i}\left(\boldsymbol{\theta}_{i}\left(t\right)\right)\right]-\widetilde{f}_{i}\left(\widetilde{\boldsymbol{\theta}}^{*}\right)\right) & \leq\frac{2MR^{2}}{\alpha T}+M\frac{\delta^{2}}{\alpha}+MC^{2}\frac{\alpha}{\beta^{2}}+MC\frac{\delta}{\beta}.\label{eq:genral_error-2}
\end{align}
\end{lem}
The proof of Lemma~\ref{lem:general_lem2} is presented in Appendix~\ref{subsec:appx_h}.

From (\ref{eq:genral_error-1}), (\ref{eq:genral_error-2}), and consider
$\delta=C\left(2+N\rho\frac{\eta}{1-\eta}\right)\frac{\alpha}{\beta}$,
we can obtain an upper bound of $\mathbb{E}[f(\overline{\boldsymbol{\theta}}^{\left(T\right)})-f(\boldsymbol{\theta}^{*})]$
as the function of the step-sizes $\alpha$ and $\beta$, i.e.,

\begin{align}
 & \mathbb{E}\left[f\left(\overline{\boldsymbol{\theta}}^{\left(T\right)}\right)-f\left(\boldsymbol{\theta}^{*}\right)\right]\nonumber \\
 & \leq\frac{2MR^{2}}{\alpha T}+MC^{2}\left(7+\frac{5N\rho\eta}{1-\eta}+\left(\frac{N\rho\eta}{1-\eta}\right)^{2}\right)\frac{\alpha}{\beta^{2}}+2\ell\beta+C\ell\left(2+\frac{N\rho\eta}{1-\eta}\right)\frac{\alpha}{\beta}\label{eq:general_upper_}\\
 & \propto O\left(\frac{1}{\alpha T}\right)+O\left(\frac{\alpha}{\beta^{2}}\right)+O\left(\beta\right)+O\left(\frac{\alpha}{\beta}\right)\label{eq:general_upper}
\end{align}

The final task is to find the reasonable values of $\alpha$ and $\beta$,
to make the above upper bound as tight as possible. Intuitively, the
desirable values of $\alpha$ and $\beta$ should be small and can
be seen as some function of $T$. In this situation $\alpha/\beta^{2}$
should be much larger than $\alpha/\beta$, which means that the upper
bound (\ref{eq:general_upper}) is not dominated by $O\left(\alpha/\beta\right)$.
For this reason, we focus on the minimization of the remaining term
$O\left(1/(\alpha T)\right)+O\left(\alpha/\beta^{2}\right)+O\left(\beta\right)$.
We have the following achievable lower bound 
\begin{align}
 & 2MR^{2}\frac{1}{\alpha T}+MC^{2}\left(7+5N\rho\frac{\eta}{1-\eta}+\left(N\rho\frac{\eta}{1-\eta}\right)^{2}\right)\frac{\alpha}{\beta^{2}}+2\ell\beta\nonumber \\
 & \geq2MRC\sqrt{14+10N\rho\frac{\eta}{1-\eta}+2\left(N\rho\frac{\eta}{1-\eta}\right)^{2}}\frac{1}{\sqrt{T}\beta}+2\ell\beta\label{eq:alpha}\\
 & \geq4\sqrt{MRC\ell\sqrt{14+10N\rho\frac{\eta}{1-\eta}+2\left(N\rho\frac{\eta}{1-\eta}\right)^{2}}}T^{-\frac{1}{4}}\propto O\left(\sqrt{M}T^{-\frac{1}{4}}\right),\label{eq:beta}
\end{align}
where both (\ref{eq:alpha}) and (\ref{eq:beta}) come from the fact
that $ax+bx^{-1}\geq2\sqrt{ab}$, $\forall a,b,x\in\mathbb{R}^{+}$
and the equality holds if and only if $ax=bx^{-1}$. Thus the equality
of (\ref{eq:alpha}) and (\ref{eq:beta}) hold only if 
\begin{equation}
\begin{cases}
2MR^{2}\frac{1}{\alpha^{*}T}=MC^{2}\left(7+5N\rho\frac{\eta}{1-\eta}+\left(N\rho\frac{\eta}{1-\eta}\right)^{2}\right)\frac{\alpha^{*}}{(\beta^{*})^{2}}\\
2MRC\sqrt{14+10N\rho\frac{\eta}{1-\eta}+2\left(N\rho\frac{\eta}{1-\eta}\right)^{2}}\frac{1}{\sqrt{T}\beta^{*}}=2\ell\beta^{*}
\end{cases}
\end{equation}
which can be solved to get $\alpha^{*}\propto T^{-\frac{3}{4}}$ and
$\beta^{*}\propto T^{-\frac{1}{4}}$ with exact expression given in
(\ref{eq:ab_opt}). With $\alpha^{*}$ and $\beta^{*}$, we can see
that the last term $O\left(\alpha/\beta\right)=O(T^{-\frac{1}{2}})$,
which is indeed much smaller than $O(T^{-\frac{1}{4}})$ when $T$
is large. In fact
\begin{align}
C\ell\left(2+N\rho\frac{\eta}{1-\eta}\right)\frac{\alpha^{*}}{\beta^{*}} & =R\ell\sqrt{\frac{2\left(2+N\rho\frac{\eta}{1-\eta}\right)^{2}}{\left(1+\left(2+N\rho\frac{\eta}{1-\eta}\right)+\left(2+N\rho\frac{\eta}{1-\eta}\right)^{2}\right)}}T^{-\frac{1}{2}}\nonumber \\
 & \leq\sqrt{2}R\ell T^{-\frac{1}{2}}.\label{eq:last_term}
\end{align}
Introduce (\ref{eq:beta}) and (\ref{eq:last_term}) into (\ref{eq:general_upper_}),
we get $\mathbb{E}[f(\overline{\boldsymbol{\theta}}^{\left(T\right)})-f\left(\boldsymbol{\theta}^{*}\right)]\propto O\left(\sqrt{M}T^{-\frac{1}{4}}\right)+O\left(T^{-\frac{1}{2}}\right)$,
with closed form expression given in (\ref{eq:general_result}).

\subsection{Discussion}

In this section, we have investigated the convergence rate of the
proposed algorithm for general Lipschitz convex function. Our result
can be seen as an extension of the classical work \cite{flaxman2005online},
which has shown that the optimization error is $O(\sqrt{M}T^{-\frac{1}{4}})$
in a centralized setting.

\section{Conclusion\label{sec:Conclusion}}

In this work, we have addressed a distributed optimization problem
in a time-varying network with non always strongly connectivity. Unlike
the classical problems where gradient information is directly available,
we consider a zeroth-order stochastic convex optimization setting.
Each node can only use a numerical observation of its local cost function
to get a biased estimation of gradient. We take a first step towards
this challenging problem. A simple distributed algorithm is considered,
with the best attainable convergence rate of $O\left(T^{-\frac{1}{2}}\right)$
after $T$ iterations, under the assumption that the global cost function
is strongly convex and local cost functions are smooth. Note that
such convergence rate is optimal in the ZOSCO problem. 

%\bibliographystyle{plain}
%\bibliography{BiblioWenjie}

\appendix
%dummy comment inserted by tex2lyx to ensure that this paragraph is not empty%dummy comment inserted by tex2lyx to ensure that this paragraph is not empty%dummy comment inserted by tex2lyx to ensure that this paragraph is not empty

\section{\label{subsec:appx_a}Proof of Lemma~\ref{lem:bias_bound}}

\label{prof_bias}We first evaluate the expected value of $\widehat{\boldsymbol{g}}_{i}\left(t\right)$
for a given $\boldsymbol{\theta}_{i}\left(t\right)$, in order to
find out its difference to the gradient $\nabla f_{i}\left(\boldsymbol{\theta}_{i}\left(t\right)\right)$.
We have 
\begin{align}
 & \mathbb{E}\left[\widehat{\boldsymbol{g}}_{i}\left(t\right)\mid\boldsymbol{\theta}_{i}\left(t\right)\right]\nonumber \\
 & =\mathbb{E}_{\boldsymbol{\nu}_{i,t},\boldsymbol{\xi}_{i,t}}\left[\beta_{t}^{-1}\boldsymbol{\nu}_{i,t}F_{i}\left(\boldsymbol{\theta}_{i}\left(t\right)+\beta_{t}\boldsymbol{\nu}_{i,t};\boldsymbol{\xi}_{i,t}\right)\mid\boldsymbol{\theta}_{i}\left(t\right)\right]\nonumber \\
 & =\beta_{t}^{-1}\mathbb{E}_{\boldsymbol{\nu}_{i,t}}\left[\boldsymbol{\nu}_{i,t}f_{i}\left(\boldsymbol{\theta}_{i}\left(t\right)+\beta_{t}\boldsymbol{\nu}_{i,t}\right)\mid\boldsymbol{\theta}_{i}\left(t\right)\right]\label{eq:rp_1}\\
 & =\beta_{t}^{-1}\mathbb{E}_{\boldsymbol{\nu}_{i,t}}\!\!\left[\boldsymbol{\nu}_{i,t}\!\!\left(\!\!f_{i}(\boldsymbol{\theta}_{i}\!\left(t\right))\!+\!\beta_{t}\boldsymbol{\nu}_{i,t}^{\mathsf{Tr}}\nabla\!f_{i}(\boldsymbol{\theta}_{i}\!\left(t\right))\!+\!\frac{1}{2}\beta_{t}^{2}\boldsymbol{\nu}_{i,t}^{\mathsf{Tr}}\nabla^{2}\!f_{i}(\widetilde{\boldsymbol{\theta}}_{i}\!\left(t\right))\boldsymbol{\nu}_{i,t}\!\!\right)\!\!\mid\!\boldsymbol{\theta}_{i}\!\left(t\right)\right]\label{eq:rp_2}\\
 & =\beta_{t}^{-1}\mathbb{E}\left[\boldsymbol{\nu}_{i,t}\right]\nabla f_{i}\left(\boldsymbol{\theta}_{i}\left(t\right)\right)+\mathbb{E}\left[\boldsymbol{\nu}_{i,t}\boldsymbol{\nu}_{i,t}^{\mathsf{Tr}}\right]\nabla f_{i}\left(\boldsymbol{\theta}_{i}\left(t\right)\right)+\boldsymbol{b}_{i}\left(t\right)\label{eq:rp_3}\\
 & =\nabla f_{i}\left(\boldsymbol{\theta}_{i}\left(t\right)\right)+\boldsymbol{b}_{i}\left(t\right).\label{eq:rp_4}
\end{align}
Note that (\ref{eq:rp_1}) is by $f_{i}\left(\boldsymbol{\theta}\right)=\mathbb{E}_{\boldsymbol{\xi}_{i,t}}\left[F_{i}\left(\boldsymbol{\theta};\boldsymbol{\xi}_{i,t}\right)\right]$
for any $\boldsymbol{\theta}\in\mathcal{K}$. (\ref{eq:rp_2}) can
be obtained by applying Taylor's theorem and mean-valued theorem,
where $\widetilde{\boldsymbol{\theta}}_{i}\left(t\right)$ locates
between $\boldsymbol{\theta}_{i}\left(t\right)$ and $\boldsymbol{\theta}_{i}\left(t\right)+\beta_{t}\boldsymbol{\nu}_{i,t}$.
In (\ref{eq:rp_3}), we introduce
\begin{equation}
\boldsymbol{b}_{i}\left(t\right)=\frac{1}{2}\beta_{t}\mathbb{E}_{\boldsymbol{\nu}_{i,t}}\left[\boldsymbol{\nu}_{i,t}\left(\boldsymbol{\nu}_{i,t}^{\mathsf{Tr}}\nabla^{2}f_{i}\left(\widetilde{\boldsymbol{\theta}}_{i}\left(t\right)\right)\boldsymbol{\nu}_{i,t}\right)\mid\boldsymbol{\theta}_{i}\left(t\right)\right].\label{eq:b_def}
\end{equation}
(\ref{eq:rp_4}) comes from the statistical property of $\boldsymbol{\nu}_{i,t}$,
i.e., $\mathbb{E}\left[\boldsymbol{\nu}_{i,t}\right]=\mathbf{0}$
and $\mathbb{E}\left[\boldsymbol{\nu}_{i,t}\boldsymbol{\nu}_{i,t}^{\mathsf{Tr}}\right]=\mathbf{I}_{M}$.
One can clearly see that $\boldsymbol{b}_{i}\left(t\right)$ defined
in (\ref{eq:b_def}) is exactly the difference between $\mathbb{E}\left[\widehat{\boldsymbol{g}}_{i}\left(t\right)\mid\boldsymbol{\theta}_{i}\left(t\right)\right]$
and $\nabla f_{i}\left(\boldsymbol{\theta}_{i}\left(t\right)\right)$.
Thus $\boldsymbol{b}_{i}\left(t\right)$ can be named as the estimation
bias of the gradient estimator. 

Our next target is to get an upper bound of $\left\Vert \boldsymbol{b}_{i}\left(t\right)\right\Vert $.
By the convexity and smoothness of $f_{i}$, we have 
\begin{align}
f_{i}\left(\boldsymbol{\theta}_{i}\left(t\right)+\beta_{t}\boldsymbol{\nu}_{i,t}\right) & \geq f_{i}\left(\boldsymbol{\theta}_{i}\left(t\right)\right)+\beta_{t}\boldsymbol{\nu}_{i,t}^{\mathsf{Tr}}\nabla f_{i}\left(\boldsymbol{\theta}_{i}\left(t\right)\right);\label{eq:concave-f}\\
f_{i}\left(\boldsymbol{\theta}_{i}\left(t\right)+\beta_{t}\boldsymbol{\nu}_{i,t}\right) & \leq f_{i}\left(\boldsymbol{\theta}_{i}\left(t\right)\right)+\beta_{t}\boldsymbol{\nu}_{i,t}^{\mathsf{Tr}}\nabla f_{i}\left(\boldsymbol{\theta}_{i}\left(t\right)\right)+\frac{1}{2}L_{i}\left\Vert \beta_{t}\boldsymbol{\nu}_{i,t}\right\Vert ^{2}.\label{eq:smooth_a}
\end{align}
By comparing (\ref{eq:concave-f})-(\ref{eq:smooth_a}) with (\ref{eq:rp_2}),
we get $0\leq\frac{1}{2}\beta_{t}^{2}\boldsymbol{\nu}_{i,t}^{\mathsf{Tr}}\nabla^{2}f_{i}(\widetilde{\boldsymbol{\theta}}_{i}\left(t\right))\boldsymbol{\nu}_{i,t}\leq\frac{1}{2}L_{i}\left\Vert \beta_{t}\boldsymbol{\nu}_{i,t}\right\Vert ^{2}$,
which means that 
\begin{equation}
\left|\boldsymbol{\nu}_{i,t}^{\mathsf{Tr}}\nabla^{2}f_{i}\left(\widetilde{\boldsymbol{\theta}}_{i}\left(t\right)\right)\boldsymbol{\nu}_{i,t}\right|\leq L_{i}\left\Vert \boldsymbol{\nu}_{i,t}\right\Vert ^{2}\leq MV^{2}L_{i}.\label{eq:bound_second_o}
\end{equation}
Thus 
\begin{align}
\left\Vert \boldsymbol{b}_{i}\left(t\right)\right\Vert  & =\frac{1}{2}\beta_{t}\left\Vert \mathbb{E}_{\boldsymbol{\nu}_{i,t}}\left[\boldsymbol{\nu}_{i,t}\left(\boldsymbol{\nu}_{i,t}^{\mathsf{Tr}}\nabla^{2}f_{i}\left(\widetilde{\boldsymbol{\theta}}_{i}\left(t\right)\right)\boldsymbol{\nu}_{i,t}\right)\mid\boldsymbol{\theta}_{i}\left(t\right)\right]\right\Vert \nonumber \\
 & \leq\frac{1}{2}\beta_{t}\sqrt{\mathbb{E}_{\boldsymbol{\nu}_{i,t}}\left[\left\Vert \boldsymbol{\nu}_{i,t}\right\Vert ^{2}\left|\boldsymbol{\nu}_{i,t}^{\mathsf{Tr}}\nabla^{2}f_{i}\left(\widetilde{\boldsymbol{\theta}}_{i}\left(t\right)\right)\boldsymbol{\nu}_{i,t}\right|^{2}\mid\boldsymbol{\theta}_{i}\left(t\right)\right]}\label{eq:rp_5}\\
 & \leq\frac{1}{2}\beta_{t}\sqrt{\mathbb{E}_{\boldsymbol{\nu}_{i,t}}\left[\left\Vert \boldsymbol{\nu}_{i,t}\right\Vert ^{2}M^{2}V^{4}L_{i}^{2}\mid\boldsymbol{\theta}_{i}\left(t\right)\right]}=\frac{1}{2}M^{\frac{3}{2}}V^{2}L_{i}\beta_{t}\label{eq:rp_6}
\end{align}
where (\ref{eq:rp_5}) comes from the fact that $\left\Vert \mathbb{E}\left[\boldsymbol{x}\right]\right\Vert \leq\sqrt{\mathbb{E}\left[\left\Vert \boldsymbol{x}\right\Vert ^{2}\right]}$
for any random vector $\boldsymbol{x}\in\mathbb{R}^{M}$, (\ref{eq:rp_6})
is by using (\ref{eq:bound_second_o}) and $\mathbb{E}\left[\boldsymbol{\nu}_{i,t}\boldsymbol{\nu}_{i,t}^{\mathsf{Tr}}\right]=\mathbf{I}_{M}$
which means that $\mathbb{E}\left[\left\Vert \boldsymbol{\nu}_{i,t}\right\Vert ^{2}\right]=M$.

\section{\label{subsec:appx_b}Proof of Lemma~\ref{lem:dist_aver}}

Define $\boldsymbol{\Phi}\left(t,s\right)=\mathbf{A}\left(t\right)\cdot\mathbf{A}\left(t-1\right)\cdot\ldots\cdot\mathbf{A}\left(s\right)$.
Introduce 
\[
\boldsymbol{\Delta}_{i}\left(t\right)=\sum_{j=1}^{N}A_{i,j}\left(t\right)\boldsymbol{\theta}_{j}\left(t\right)-\underbrace{\Pi_{\mathcal{K}_{t+1}}\left(\sum_{j=1}^{N}A_{i,j}\left(t\right)\boldsymbol{\theta}_{j}\left(t\right)-\alpha_{t}\widehat{\boldsymbol{g}}_{i}\left(t\right)\right)}_{=\boldsymbol{\theta}_{i}\left(t+1\right)}.
\]
Since $\boldsymbol{\theta}_{j}\left(t\right)\in\mathcal{K}_{t}\subseteq\mathcal{K}_{t+1}$,
$\forall j$, we have $\sum_{j=1}^{N}A_{i,j}\left(t\right)\boldsymbol{\theta}_{j}\left(t\right)\in\mathcal{K}_{t+1}$
as $\mathcal{K}_{t+1}$ is convex. On the other hand, we also have
$\Pi_{\mathcal{K}_{t+1}}\left(\boldsymbol{x}\right)\in\mathcal{K}_{t+1}$.
Thus 
\begin{align}
\left\Vert \boldsymbol{\Delta}_{i}\left(t\right)\right\Vert  & =\left\Vert \sum_{j=1}^{N}A_{i,j}\left(t\right)\boldsymbol{\theta}_{j}\left(t\right)-\Pi_{\mathcal{K}_{t+1}}\left(\sum_{j=1}^{N}A_{i,j}\left(t\right)\boldsymbol{\theta}_{j}\left(t\right)-\alpha_{t}\widehat{\boldsymbol{g}}_{i}\left(t\right)\right)\right\Vert \nonumber \\
 & \leq\left\Vert \sum_{j=1}^{N}A_{i,j}\left(t\right)\boldsymbol{\theta}_{j}\left(t\right)-\left(\sum_{j=1}^{N}A_{i,j}\left(t\right)\boldsymbol{\theta}_{j}\left(t\right)-\alpha_{t}\widehat{\boldsymbol{g}}_{i}\left(t\right)\right)\right\Vert \nonumber \\
 & =\alpha_{t}\left\Vert \widehat{\boldsymbol{g}}_{i}\left(t\right)\right\Vert =\frac{\alpha_{t}}{\beta_{t}}\left|F_{i}\left(\boldsymbol{\theta}_{i}\left(t\right)+\beta_{t}\boldsymbol{\Phi}_{i}\left(t\right);\boldsymbol{\xi}_{i,t}\right)\right|\left\Vert \boldsymbol{\nu}_{i,t}\right\Vert \nonumber \\
 & \leq\frac{\alpha_{t}}{\beta_{t}}C\left\Vert \boldsymbol{\nu}_{i,t}\right\Vert .\label{eq:bound_delta}
\end{align}
recall that $\left|F_{i}\right|\leq C$ by Assumption~A3. From the
definition of $\boldsymbol{\Delta}_{i}\left(t\right)$, we evaluate
\begin{align*}
 & \boldsymbol{\theta}_{i}\left(t\right)=\sum_{j_{1}=1}^{N}A_{i,j_{1}}\left(t-1\right)\boldsymbol{\theta}_{j_{1}}\left(t-1\right)-\boldsymbol{\Delta}_{i}\left(t-1\right)\\
 & =\sum_{j_{1}=1}^{N}A_{i,j_{1}}\left(t-1\right)\left(\sum_{j_{2}=1}^{N}A_{j_{1},j_{2}}\left(t-2\right)\boldsymbol{\theta}_{j_{2}}\left(t-2\right)-\boldsymbol{\Delta}_{j_{1}}\left(t-2\right)\right)-\boldsymbol{\Delta}_{i}\left(t-1\right)=\cdots\\
 & =\!\sum_{j_{t}=1}^{N}\!\left[\boldsymbol{\Phi}\left(t-1,0\right)\right]_{i,j_{t}}\!\boldsymbol{\theta}_{j_{t}}\!\left(0\right)\!-\!\sum_{s=1}^{t-1}\sum_{j_{s}=1}^{N}\!\left[\boldsymbol{\Phi}\left(t-1,t-s\right)\right]_{i,j_{s}}\!\boldsymbol{\Delta}_{j_{s}}\!\left(t-s-1\right)\!-\!\boldsymbol{\Delta}_{i}\!\left(t-1\right)\\
 & =\sum_{j=1}^{N}\left[\boldsymbol{\Phi}\left(t-1,0\right)\right]_{i,j}\boldsymbol{\theta}_{j}\left(0\right)-\sum_{s=1}^{t-1}\sum_{j=1}^{N}\left[\boldsymbol{\Phi}\left(t-1,s\right)\right]_{i,j}\boldsymbol{\Delta}_{j}\left(s-1\right)-\boldsymbol{\Delta}_{i}\left(t-1\right)
\end{align*}
Meanwhile 
\begin{align*}
 & \overline{\boldsymbol{\theta}}\left(t\right)=\frac{1}{N}\sum_{i=1}^{N}\boldsymbol{\theta}_{i}\left(t\right)=\frac{1}{N}\sum_{i=1}^{N}\left(\sum_{j=1}^{N}A_{i,j}\left(t-1\right)\boldsymbol{\theta}_{j}\left(t-1\right)-\boldsymbol{\Delta}_{i}\left(t-1\right)\right)\\
 & =\frac{1}{N}\sum_{j=1}^{N}\!\left(\sum_{i=1}^{N}\!A_{i,j}\!\left(t-1\right)\!\right)\!\boldsymbol{\theta}_{j}\!\left(t-1\right)\!-\!\frac{1}{N}\sum_{i=1}^{N}\!\boldsymbol{\Delta}_{i}\!\left(t-1\right)=\overline{\boldsymbol{\theta}}\!\left(t-1\right)-\!\frac{1}{N}\sum_{i=1}^{N}\!\boldsymbol{\Delta}_{i}\!\left(t-1\right)\\
 & =\ldots=\frac{1}{N}\sum_{j=1}^{N}\boldsymbol{\theta}_{j}\left(0\right)-\sum_{s=1}^{t-1}\sum_{j=1}^{N}\frac{1}{N}\boldsymbol{\Delta}_{j}\left(s-1\right)-\frac{1}{N}\sum_{j=1}^{N}\boldsymbol{\Delta}_{j}\left(t-1\right)
\end{align*}
Now we can have the difference between $\boldsymbol{\theta}_{i}\left(t\right)$
and $\overline{\boldsymbol{\theta}}\left(t\right)$, since $\boldsymbol{\theta}_{j}\left(0\right)=\mathbf{0}$,
$\forall j$:
\begin{align}
\left\Vert \boldsymbol{\theta}_{i}\left(t\right)-\overline{\boldsymbol{\theta}}\left(t\right)\right\Vert  & \leq\left\Vert -\boldsymbol{\Delta}_{i}\left(t-1\right)+\frac{1}{N}\sum_{j=1}^{N}\boldsymbol{\Delta}_{j}\left(t-1\right)\right\Vert \nonumber \\
 & \qquad+\sum_{s=1}^{t-1}\sum_{j=1}^{N}\left|\left[\boldsymbol{\Phi}\left(t-1,s\right)\right]_{i,j}-\frac{1}{N}\right|\left\Vert \boldsymbol{\Delta}_{j}\left(s-1\right)\right\Vert \label{eq:xi_diff}
\end{align}
Under Assumption A1, according Proposition~1 in the reference \cite{nedic2009distributed},
we have
\[
\left|\left[\boldsymbol{\Phi}\left(t,s\right)\right]_{i,j}-\frac{1}{N}\right|\leq\rho\eta^{t-s}
\]
where $\rho=2\frac{1+a^{\left(N-1\right)\tau}}{1-a^{\left(N-1\right)\tau}}$
and $\eta=\left(1-a^{\left(N-1\right)\tau}\right)^{\frac{1}{\left(N-1\right)\tau}}$
as introduces in Lemma~\ref{lem:dist_aver}.

Consider the first case where $\left\Vert \boldsymbol{\nu}_{i,t}\right\Vert \leq\sqrt{M}V$
and the step-sizes are vanishing. From (\ref{eq:xi_diff}) and (\ref{eq:bound_delta}),
we have 
\begin{align}
 & \left\Vert \boldsymbol{\theta}_{i}\left(t\right)-\overline{\boldsymbol{\theta}}\left(t\right)\right\Vert \nonumber \\
 & \leq\left\Vert \boldsymbol{\Delta}_{i}\left(t-1\right)\right\Vert +\frac{1}{N}\sum_{j=1}^{N}\left\Vert \boldsymbol{\Delta}_{j}\left(t-1\right)\right\Vert +\sum_{s=1}^{t-1}\sum_{j=1}^{N}\rho\eta^{t-s-1}\left\Vert \boldsymbol{\Delta}_{j}\left(s-1\right)\right\Vert \label{eq:xi_diff-1}\\
 & \leq2C\left\Vert \boldsymbol{\nu}_{i,t}\right\Vert \frac{\alpha_{t-1}}{\beta_{t-1}}+N\rho C\left\Vert \boldsymbol{\nu}_{i,t}\right\Vert \sum_{s=0}^{t-2}\frac{\alpha_{s}}{\beta_{s}}\eta^{t-s-2}\nonumber \\
 & \leq\sqrt{M}VC\left(2^{1+c_{3}}+N\rho\left(\frac{1+c_{3}}{1-\eta}+\frac{c_{3}}{\left(1-\eta\right)^{2}}\right)\right)\frac{\alpha_{t}}{\beta_{t}}\label{eq:theta_diff}
\end{align}
Denote $\gamma_{s}=\alpha_{s}/\beta_{s}=\gamma_{0}s^{-c_{3}}$, then
(\ref{eq:theta_diff}) comes from $\frac{\gamma_{t-1}}{\gamma_{t}}\leq2^{c_{3}}$
and $\sum_{s=0}^{t-2}\frac{\alpha_{s}}{\beta_{s}}\eta^{t-s-2}\leq\gamma_{t}\left(\frac{1+c_{3}}{1-\eta}+\frac{c_{3}}{\left(1-\eta\right)^{2}}\right)$
which can be proved by the following:
\begin{align}
 & \sum_{s=0}^{t-2}\eta^{t-s-2}\frac{\alpha_{s}}{\beta_{s}}=\gamma_{t}\sum_{s=0}^{t-2}\eta^{t-s-2}\frac{\gamma_{s}}{\gamma_{t}}=\gamma_{t}\sum_{s=0}^{t-2}\eta^{t-s-2}\left(\frac{s}{t}\right)^{-c_{3}}\nonumber \\
 & \leq\gamma_{t}\left(\eta^{t-0-2}+\sum_{s=1}^{t-2}\eta^{t-s-2}\left(1+c_{3}\frac{t-s}{s}\right)\right)\label{eq:tt1}\\
 & \leq\gamma_{t}\left(\sum_{s=0}^{t-2}\eta^{t-s-2}+c_{3}\sum_{s=1}^{t-2}\eta^{t-s-2}\left(t-s\right)\right)\nonumber \\
 & =\gamma_{t}\left(\sum_{s=0}^{t-2}\eta^{s}+\frac{c_{3}}{\eta}\sum_{s=2}^{t-1}\eta^{s-1}s\right)=\gamma_{t}\left(\sum_{s=0}^{t-2}\eta^{s}+\frac{c_{3}}{\eta}\frac{d}{d\eta}\left(\sum_{s=2}^{t-1}\eta^{s}\right)\right)\nonumber \\
 & =\gamma_{t}\left(\frac{1-\eta^{t-1}}{1-\eta}+\frac{c_{3}}{\eta}\frac{\left(1-\eta\right)\left(2\eta-t\eta^{t-1}\right)+\left(\eta^{2}-\eta^{t}\right)}{\left(1-\eta\right)^{2}}\right)\nonumber \\
 & \leq\gamma_{t}\left(\frac{1}{1-\eta}+c_{3}\frac{2\left(1-\eta\right)+\eta}{\left(1-\eta\right)^{2}}\right)=\gamma_{t}\left(\frac{1+c_{3}}{1-\eta}+\frac{c_{3}}{\left(1-\eta\right)^{2}}\right),\nonumber 
\end{align}
note that (\ref{eq:tt1}) is by $\left(1+x\right)^{c_{3}}\leq1+c_{3}x$
for any $x>0$ and $c_{3}<1$.

In the end, we consider the more simple case where $\left\Vert \boldsymbol{\nu}_{i,t}\right\Vert =1$
and the step-sizes are constant. With similar steps as (\ref{eq:theta_diff}),
we have 
\begin{align}
\left\Vert \boldsymbol{\theta}_{i}\left(t\right)-\overline{\boldsymbol{\theta}}\left(t\right)\right\Vert  & \leq\frac{\alpha}{\beta}C\left(2+N\rho\sum_{s=0}^{t-2}\eta^{t-s-2}\right)=\frac{\alpha}{\beta}C\left(2+N\rho\eta\frac{1-\eta^{t-1}}{1-\eta}\right)\nonumber \\
 & \leq\frac{\alpha}{\beta}C\left(2+N\rho\frac{\eta}{1-\eta}\right)\label{eq:theta_diff-1}
\end{align}
which concludes the proof.

\section{\label{subsec:appx_c}Proof of inequality (\ref{eq:error_divergence})}

We have 
\begin{align}
 & f\left(\frac{1}{NT}\sum_{t=1}^{T}\sum_{i=1}^{N}\boldsymbol{\theta}_{i}\left(t\right)\right)-f\left(\boldsymbol{\theta}^{*}\right)\leq\frac{1}{T}\sum_{t=1}^{T}f\left(\frac{1}{N}\sum_{i=1}^{N}\boldsymbol{\theta}_{i}\left(t\right)\right)-f\left(\boldsymbol{\theta}^{*}\right)\label{eq:p1}\\
 & =\frac{1}{T}\sum_{t=1}^{T}\left(f\left(\overline{\boldsymbol{\theta}}\left(t\right)\right)-f\left(\boldsymbol{\theta}^{*}\right)\right)=\frac{1}{T}\sum_{t=1}^{T}\frac{1}{N}\sum_{i=1}^{N}\left(f_{i}\left(\overline{\boldsymbol{\theta}}\left(t\right)\right)-f_{i}\left(\boldsymbol{\theta}^{*}\right)\right)\nonumber \\
 & \leq\frac{1}{T}\sum_{t=1}^{T}\frac{1}{N}\sum_{i=1}^{N}\frac{L_{i}}{2}\left\Vert \overline{\boldsymbol{\theta}}\left(t\right)-\boldsymbol{\theta}^{*}\right\Vert ^{2}=\frac{1}{T}\sum_{t=1}^{T}\frac{L}{2}\left\Vert \overline{\boldsymbol{\theta}}\left(t\right)-\boldsymbol{\theta}^{*}\right\Vert ^{2}\label{eq:p2}\\
 & \leq\frac{L}{2T}\sum_{t=1}^{T}\frac{1}{N}\sum_{i=1}^{N}\left\Vert \boldsymbol{\theta}_{i}\left(t\right)-\boldsymbol{\theta}^{*}\right\Vert ^{2},\label{eq:p3}
\end{align}
where (\ref{eq:p1}) is by the convexity of $f$; (\ref{eq:p2}) is
by the assumption that $f_{i}$ is $L_{i}-$smooth and $L=\frac{1}{N}\sum_{i=1}^{N}L_{i}$;
(\ref{eq:p3}) comes from
\begin{equation}
\left\Vert \overline{\boldsymbol{\theta}}\left(t\right)-\boldsymbol{\theta}^{*}\right\Vert ^{2}=\left\Vert \frac{1}{N}\sum_{i=1}^{N}\left(\boldsymbol{\theta}_{i}\left(t\right)-\boldsymbol{\theta}^{*}\right)\right\Vert ^{2}\leq\frac{1}{N}\sum_{i=1}^{N}\left\Vert \boldsymbol{\theta}_{i}\left(t\right)-\boldsymbol{\theta}^{*}\right\Vert ^{2}.\label{eq:p4}
\end{equation}
By taking expectation on both sides of the inequality, we get (\ref{eq:error_divergence}).

\section{\label{subsec:appx_d}Proof of Lemma~\ref{lem:algo_iter_strong}}

\label{proof_iter}We need to evaluate
\begin{align}
 & \sum_{i=1}^{N}\left\Vert \boldsymbol{\theta}_{i}\left(t+1\right)-\boldsymbol{\theta}^{*}\right\Vert ^{2}=\sum_{i=1}^{N}\left\Vert \Pi_{\mathcal{K}_{t+1}}\left(\sum_{j=1}^{N}A_{i,j}\left(t\right)\boldsymbol{\theta}_{j}\left(t\right)-\alpha_{t}\widehat{\boldsymbol{g}}_{i}\left(t\right)\right)-\boldsymbol{\theta}^{*}\right\Vert ^{2}\nonumber \\
 & \leq\sum_{i=1}^{N}\left\Vert \sum_{j=1}^{N}A_{i,j}\left(t\right)\boldsymbol{\theta}_{j}\left(t\right)-\alpha_{t}\widehat{\boldsymbol{g}}_{i}\left(t\right)-\boldsymbol{\theta}^{*}\right\Vert ^{2}\label{eq:tt3}\\
 & =\sum_{i=1}^{N}\left\Vert \left(\overline{\boldsymbol{\theta}}\left(t\right)-\boldsymbol{\theta}^{*}\right)+\left(\sum_{j=1}^{N}A_{i,j}\left(t\right)\boldsymbol{\theta}_{j}\left(t\right)-\alpha_{t}\widehat{\boldsymbol{g}}_{i}\left(t\right)-\overline{\boldsymbol{\theta}}\left(t\right)\right)\right\Vert ^{2}\nonumber \\
 & =N\left\Vert \overline{\boldsymbol{\theta}}\left(t\right)-\boldsymbol{\theta}^{*}\right\Vert ^{2}+\sum_{i=1}^{N}\left\Vert \sum_{j=1}^{N}A_{i,j}\left(t\right)\boldsymbol{\theta}_{j}\left(t\right)-\overline{\boldsymbol{\theta}}\left(t\right)-\alpha_{t}\widehat{\boldsymbol{g}}_{i}\left(t\right)\right\Vert ^{2}\nonumber \\
 & -2\alpha_{t}\sum_{i=1}^{N}\left\langle \overline{\boldsymbol{\theta}}\left(t\right)-\boldsymbol{\theta}^{*},\widehat{\boldsymbol{g}}_{i}\left(t\right)\right\rangle +2\sum_{i=1}^{N}\left\langle \overline{\boldsymbol{\theta}}\left(t\right)-\boldsymbol{\theta}^{*},\sum_{j=1}^{N}A_{i,j}\left(t\right)\boldsymbol{\theta}_{j}\left(t\right)-\overline{\boldsymbol{\theta}}\left(t\right)\right\rangle ,\label{eq:dt1-dt}
\end{align}
where (\ref{eq:tt3}) holds as $\boldsymbol{\theta}^{*}\in\mathcal{K}_{0}\subseteq\mathcal{K}_{t+1}$
and $\Pi_{\mathcal{K}_{t+1}}\left(\sum_{j=1}^{N}A_{i,j}\left(t\right)\boldsymbol{\theta}_{j}\left(t\right)-\alpha_{t}\widehat{\boldsymbol{g}}_{i}\left(t\right)\right)\subseteq\mathcal{K}_{t+1}$,
thus 
\[
\left\Vert \sum_{j=1}^{N}A_{i,j}\left(t\right)\boldsymbol{\theta}_{j}\left(t\right)-\alpha_{t}\widehat{\boldsymbol{g}}_{i}\left(t\right)-\boldsymbol{\theta}^{*}\right\Vert \geq\left\Vert \Pi_{\mathcal{K}_{t+1}}\left(\sum_{j=1}^{N}A_{i,j}\left(t\right)\boldsymbol{\theta}_{j}\left(t\right)-\alpha_{t}\widehat{\boldsymbol{g}}_{i}\left(t\right)\right)-\boldsymbol{\theta}^{*}\right\Vert .
\]

We then need to investigate the term $\sum_{i=1}^{N}\left\langle \boldsymbol{\theta}_{i}\left(t\right)-\boldsymbol{\theta}^{*},\widehat{\boldsymbol{g}}_{i}\left(t\right)\right\rangle $.
For any $i\in\mathcal{V},$ we evaluate 
\begin{align*}
\widehat{\boldsymbol{g}}_{i}\left(t\right) & =\mathbb{E}\left[\widehat{\boldsymbol{g}}_{i}\left(t\right)\mid\boldsymbol{\theta}_{i}\left(t\right)\right]+\widehat{\boldsymbol{g}}_{i}\left(t\right)-\mathbb{E}\left[\widehat{\boldsymbol{g}}_{i}\left(t\right)\mid\boldsymbol{\theta}_{i}\left(t\right)\right]\\
 & =\nabla f_{i}\left(\boldsymbol{\theta}_{i}\left(t\right)\right)+\boldsymbol{b}_{i}\left(t\right)+\boldsymbol{e}_{i}\left(t\right)
\end{align*}
where we denote $\boldsymbol{e}_{i}\left(t\right)=\widehat{\boldsymbol{g}}_{i}\left(t\right)-\mathbb{E}\left[\widehat{\boldsymbol{g}}_{i}\left(t\right)\mid\boldsymbol{\theta}_{i}\left(t\right)\right]$
and recall that 
\[
\boldsymbol{b}_{i}\left(t\right)=\mathbb{E}\left[\widehat{\boldsymbol{g}}_{i}\left(t\right)\mid\boldsymbol{\theta}_{i}\left(t\right)\right]-\nabla f_{i}\left(\boldsymbol{\theta}_{i}\left(t\right)\right).
\]
Thus we have
\begin{align}
 & \left\langle \overline{\boldsymbol{\theta}}\left(t\right)-\boldsymbol{\theta}^{*},\widehat{\boldsymbol{g}}_{i}\left(t\right)\right\rangle =\left\langle \overline{\boldsymbol{\theta}}\left(t\right)-\boldsymbol{\theta}^{*},\nabla f_{i}\left(\boldsymbol{\theta}_{i}\left(t\right)\right)\right\rangle +\left\langle \overline{\boldsymbol{\theta}}\left(t\right)-\boldsymbol{\theta}^{*},\boldsymbol{b}_{i}\left(t\right)+\boldsymbol{e}_{i}\left(t\right)\right\rangle \nonumber \\
 & =\left\langle \overline{\boldsymbol{\theta}}\left(t\right)-\boldsymbol{\theta}^{*},\nabla f_{i}\left(\overline{\boldsymbol{\theta}}\left(t\right)\right)\right\rangle +\left\langle \overline{\boldsymbol{\theta}}\left(k\right)-\boldsymbol{\theta}^{*},\nabla f_{i}\left(\boldsymbol{\theta}_{i}\left(t\right)\right)-\nabla f_{i}\left(\overline{\boldsymbol{\theta}}\left(t\right)\right)\right\rangle \nonumber \\
 & \qquad+\left\langle \overline{\boldsymbol{\theta}}\left(t\right)-\boldsymbol{\theta}^{*},\boldsymbol{b}_{i}\left(t\right)+\boldsymbol{e}_{i}\left(t\right)\right\rangle \nonumber \\
 & \geq\left\langle \overline{\boldsymbol{\theta}}\left(t\right)-\boldsymbol{\theta}^{*},\nabla f_{i}\left(\overline{\boldsymbol{\theta}}\left(t\right)\right)\right\rangle +\left\langle \overline{\boldsymbol{\theta}}\left(t\right)-\boldsymbol{\theta}^{*},\boldsymbol{e}_{i}\left(t\right)\right\rangle \label{eq:xg11}\\
 & \qquad-\left\Vert \overline{\boldsymbol{\theta}}\left(t\right)-\boldsymbol{\theta}^{*}\right\Vert \left(\left\Vert \nabla f_{i}\left(\boldsymbol{\theta}_{i}\left(t\right)\right)-\nabla f_{i}\left(\overline{\boldsymbol{\theta}}\left(t\right)\right)\right\Vert +\left\Vert \boldsymbol{b}_{i}\left(t\right)\right\Vert \right)\nonumber \\
 & \geq\left\langle \overline{\boldsymbol{\theta}}\left(t\right)-\boldsymbol{\theta}^{*},\nabla f_{i}\left(\overline{\boldsymbol{\theta}}\left(t\right)\right)\right\rangle +\left\langle \overline{\boldsymbol{\theta}}\left(t\right)-\boldsymbol{\theta}^{*},\boldsymbol{e}_{i}\left(t\right)\right\rangle \label{eq:xg1}\\
 & \qquad-\left\Vert \overline{\boldsymbol{\theta}}\left(t\right)-\boldsymbol{\theta}^{*}\right\Vert \left(L_{i}\left\Vert \boldsymbol{\theta}_{i}\left(t\right)-\overline{\boldsymbol{\theta}}\left(t\right)\right\Vert +\left\Vert \boldsymbol{b}_{i}\left(t\right)\right\Vert \right)\nonumber 
\end{align}
where (\ref{eq:xg11}) is by $\left|\left\langle \boldsymbol{a}_{1},\boldsymbol{a}_{2}\right\rangle \right|\leq\left\Vert \boldsymbol{a}_{1}\right\Vert \left\Vert \boldsymbol{a}_{2}\right\Vert $
for any vectors $\boldsymbol{a}_{1}$ and $\boldsymbol{a}_{2}$ of
the same dimension; (\ref{eq:xg1}) is by the assumption that each
average local cost function $f_{i}$ is $L_{i}-$smooth. Based on
(\ref{eq:xg1}), we have 
\begin{align}
 & \sum_{i=1}^{N}\left\langle \overline{\boldsymbol{\theta}}\left(t\right)-\boldsymbol{\theta}^{*},\widehat{\boldsymbol{g}}_{i}\left(t\right)\right\rangle \geq\left\langle \overline{\boldsymbol{\theta}}\left(t\right)-\boldsymbol{\theta}^{*},\sum_{i=1}^{N}\nabla f_{i}\left(\overline{\boldsymbol{\theta}}\left(t\right)\right)\right\rangle +\sum_{i=1}^{N}\left\langle \overline{\boldsymbol{\theta}}\left(t\right)-\boldsymbol{\theta}^{*},\boldsymbol{e}_{i}\left(t\right)\right\rangle \nonumber \\
 & \qquad-\left\Vert \overline{\boldsymbol{\theta}}\left(t\right)-\boldsymbol{\theta}^{*}\right\Vert \sum_{i=1}^{N}\left(L_{i}\left\Vert \boldsymbol{\theta}_{i}\left(t\right)-\overline{\boldsymbol{\theta}}\left(t\right)\right\Vert +\left\Vert \boldsymbol{b}_{i}\left(t\right)\right\Vert \right)\nonumber \\
 & \geq N\mu\left\Vert \overline{\boldsymbol{\theta}}\left(t\right)-\boldsymbol{\theta}^{*}\right\Vert ^{2}\!-\!N\!(\frac{M^{\frac{3}{2}}}{2}V^{2}L\beta_{t}\!+\!L\delta_{t})\!\left\Vert \overline{\boldsymbol{\theta}}\left(t\right)-\boldsymbol{\theta}^{*}\right\Vert \!+\!\sum_{i=1}^{N}\!\left\langle \overline{\boldsymbol{\theta}}\left(t\right)\!-\boldsymbol{\theta}^{*},\boldsymbol{e}_{i}\left(t\right)\right\rangle \label{eq:xg2}
\end{align}
where $\sum_{i=1}^{N}\left\Vert \boldsymbol{b}_{i}\left(t\right)\right\Vert \leq\frac{1}{2}NM^{\frac{3}{2}}V^{2}L\beta_{t}$
by using Lemma~\ref{lem:bias_bound}, $\sum_{i=1}^{N}L_{i}\left\Vert \boldsymbol{\theta}_{i}\left(t\right)-\overline{\boldsymbol{\theta}}\left(t\right)\right\Vert \leq\sum_{i=1}^{N}L_{i}\delta_{t}=NL\delta_{t}$
by Lemma~\ref{lem:dist_aver}, and 
\[
\left\langle \overline{\boldsymbol{\theta}}\left(t\right)-\boldsymbol{\theta}^{*},\sum_{i=1}^{N}\nabla f_{i}\left(\overline{\boldsymbol{\theta}}\left(t\right)\right)\right\rangle =N\left\langle \overline{\boldsymbol{\theta}}\left(t\right)-\boldsymbol{\theta}^{*},\nabla f\left(\overline{\boldsymbol{\theta}}\left(t\right)\right)\right\rangle \geq N\mu\left\Vert \overline{\boldsymbol{\theta}}\left(t\right)-\boldsymbol{\theta}^{*}\right\Vert ^{2},
\]
which can be got by the assumption that $f$ is $\mu$-strongly convex
(Assumption~A4), i.e., 
\begin{equation}
\left\langle \overline{\boldsymbol{\theta}}\left(t\right)-\boldsymbol{\theta}^{*},\nabla f\left(\overline{\boldsymbol{\theta}}\left(t\right)\right)-\nabla f\left(\boldsymbol{\theta}^{*}\right)\right\rangle \geq\mu\left\Vert \overline{\boldsymbol{\theta}}\left(t\right)-\boldsymbol{\theta}^{*}\right\Vert ^{2}\label{eq:strong}
\end{equation}
and 
\begin{equation}
\left\langle \overline{\boldsymbol{\theta}}\left(t\right)-\boldsymbol{\theta}^{*},\nabla f\left(\boldsymbol{\theta}^{*}\right)\right\rangle \geq0\label{eq:mini_pro}
\end{equation}
since $\boldsymbol{\theta}^{*}$ is minimizer of the convex function
$f$.

Back to (\ref{eq:dt1-dt}), we still need to derive 
\begin{align}
 & \sum_{i=1}^{N}\!\left\langle \!\overline{\boldsymbol{\theta}}\!\left(t\right)-\boldsymbol{\theta}^{*},\sum_{j=1}^{N}\!A_{i,j}\!\left(t\right)\boldsymbol{\theta}_{j}\!\left(t\right)\!-\overline{\boldsymbol{\theta}}\!\left(t\right)\!\right\rangle \nonumber \\
 & =\left\langle \overline{\boldsymbol{\theta}}\left(t\right)-\boldsymbol{\theta}^{*},\sum_{j=1}^{N}\left(\sum_{i=1}^{N}A_{i,j}\left(t\right)\right)\boldsymbol{\theta}_{j}\left(t\right)-N\overline{\boldsymbol{\theta}}\left(t\right)\right\rangle \nonumber \\
 & =\left\langle \overline{\boldsymbol{\theta}}\left(t\right)-\boldsymbol{\theta}^{*},\sum_{j=1}^{N}\boldsymbol{\theta}_{j}\left(t\right)-N\frac{1}{N}\sum_{i=1}^{N}\boldsymbol{\theta}_{i}\left(t\right)\right\rangle =0,\label{eq:xg3}
\end{align}
thanks to the assumption that $\mathbf{A}(t)$ is doubly stochastic,
$\sum_{i=1}^{N}A_{i,j}\left(t\right)=1$. Meanwhile
\begin{align}
 & \sum_{i=1}^{N}\left\Vert \sum_{j=1}^{N}A_{i,j}\left(t\right)\boldsymbol{\theta}_{j}\left(t\right)-\overline{\boldsymbol{\theta}}\left(t\right)-\alpha_{t}\widehat{\boldsymbol{g}}_{i}\left(t\right)\right\Vert ^{2}\nonumber \\
 & \leq2\sum_{i=1}^{N}\left(\left\Vert \sum_{j=1}^{N}A_{i,j}\left(t\right)\boldsymbol{\theta}_{j}\left(t\right)-\overline{\boldsymbol{\theta}}\left(t\right)\right\Vert ^{2}+\left\Vert \alpha_{t}\widehat{\boldsymbol{g}}_{i}\left(t\right)\right\Vert ^{2}\right)\nonumber \\
 & \leq2\sum_{i=1}^{N}\left(\left(\sum_{j=1}^{N}A_{i,j}\left(t\right)\left\Vert \boldsymbol{\theta}_{j}\left(t\right)-\overline{\boldsymbol{\theta}}\left(t\right)\right\Vert \right)^{2}+\left\Vert \alpha_{t}\widehat{\boldsymbol{g}}_{i}\left(t\right)\right\Vert ^{2}\right)\nonumber \\
 & \leq2\sum_{i=1}^{N}\left(\sum_{j=1}^{N}A_{i,j}\left(t\right)\delta_{t}\right)^{2}+2\alpha_{t}^{2}\sum_{i=1}^{N}\frac{MV^{2}C^{2}}{\beta_{t}^{2}}=2N\delta_{t}^{2}+2MNV^{2}C^{2}\frac{\alpha_{t}^{2}}{\beta_{t}^{2}}\label{eq:xg4}
\end{align}

Introduce (\ref{eq:xg2}), (\ref{eq:xg3}), and (\ref{eq:xg4}) into
(\ref{eq:dt1-dt}), we get
\begin{align}
 & \sum_{i=1}^{N}\left\Vert \boldsymbol{\theta}_{i}\left(t+1\right)-\boldsymbol{\theta}^{*}\right\Vert ^{2}\nonumber \\
 & \leq N\left(1-2\mu\alpha_{t}\right)\left\Vert \overline{\boldsymbol{\theta}}\left(t\right)-\boldsymbol{\theta}^{*}\right\Vert ^{2}+NL\alpha_{t}\left(M^{\frac{3}{2}}V^{2}\beta_{t}+2\delta_{t}\right)\left\Vert \overline{\boldsymbol{\theta}}\left(t\right)-\boldsymbol{\theta}^{*}\right\Vert \nonumber \\
 & \qquad+2N\delta_{t}^{2}+2MNV^{2}C^{2}\frac{\alpha_{t}^{2}}{\beta_{t}^{2}}-2\alpha_{t}\sum_{i=1}^{N}\left\langle \overline{\boldsymbol{\theta}}\left(t\right)-\boldsymbol{\theta}^{*},\boldsymbol{e}_{i}\left(t\right)\right\rangle \nonumber \\
 & \leq\left(1-2\mu\alpha_{t}\right)\sum_{i=1}^{N}\left\Vert \boldsymbol{\theta}_{i}\left(t\right)-\boldsymbol{\theta}^{*}\right\Vert ^{2}\!+\!\sqrt{N}L\alpha_{t}(M^{\frac{3}{2}}V^{2}\beta_{t}\!+2\delta_{t})\!\!\!\sqrt{\sum_{i=1}^{N}\left\Vert \boldsymbol{\theta}_{i}\!\left(t\right)\!-\boldsymbol{\theta}^{*}\right\Vert ^{2}}\label{eq:dt1_2}\\
 & \qquad+2N\delta_{t}^{2}+2MNV^{2}C^{2}\frac{\alpha_{t}^{2}}{\beta_{t}^{2}}-2\alpha_{t}\sum_{i=1}^{N}\left\langle \overline{\boldsymbol{\theta}}\left(t\right)-\boldsymbol{\theta}^{*},\boldsymbol{e}_{i}\left(t\right)\right\rangle ,\nonumber 
\end{align}
recall that we have shown $\left\Vert \overline{\boldsymbol{\theta}}\left(t\right)-\boldsymbol{\theta}^{*}\right\Vert ^{2}\leq\frac{1}{N}\sum_{i=1}^{N}\left\Vert \boldsymbol{\theta}_{i}\left(t\right)-\boldsymbol{\theta}^{*}\right\Vert ^{2}$
in (\ref{eq:p4}).

We can take expectation on both sides of (\ref{eq:dt1_2}), considering
the the fact that $\boldsymbol{e}_{i}\left(t\right)$ has zero-mean
by definition and
\[
\mathbb{E}\left[\sqrt{\sum_{i=1}^{N}\left\Vert \boldsymbol{\theta}_{i}\left(t\right)-\boldsymbol{\theta}^{*}\right\Vert ^{2}}\right]\leq\sqrt{\mathbb{E}\left[\sum_{i=1}^{N}\left\Vert \boldsymbol{\theta}_{i}\left(t\right)-\boldsymbol{\theta}^{*}\right\Vert ^{2}\right]}=\sqrt{d_{t}}
\]
we get
\begin{equation}
d_{t+1}\leq\left(1-2\mu\alpha_{t}\right)d_{t}+\sqrt{N}L\alpha_{t}(M^{\frac{3}{2}}V^{2}\beta_{t}+2\delta_{t})\sqrt{d_{t}}+2N\delta_{t}^{2}+2MNV^{2}C^{2}\frac{\alpha_{t}^{2}}{\beta_{t}^{2}}.\label{eq:evo_strong_1}
\end{equation}
By considering the definition of $\lambda_{\mathsf{I}}$, $\lambda_{\mathsf{II}}$,
and $\lambda_{\mathsf{III}}$ in Theorem~\ref{thm:result_strong},
the inequality (\ref{eq:evo_strong_1}) can be written as (\ref{eq:evo_strong}),
which concludes the proof.

\section{\label{subsec:appx_e}Proof of Lemma~\ref{lem:rate_strong}}

This proof is by induction. Consider first the situation where $t\leq t_{0}=\left\lceil 2\mu\alpha_{0}\right\rceil $.
Since $\left\Vert \boldsymbol{\theta}\right\Vert \leq R$ for any
$\boldsymbol{\theta}\in\mathcal{K}$, we have, $\forall t\leq t_{0}$,
\begin{equation}
d_{t}=\mathbb{E}\left[\sum_{i=1}^{N}\left\Vert \boldsymbol{\theta}_{i}\left(t\right)-\boldsymbol{\theta}^{*}\right\Vert ^{2}\right]\leq\mathbb{E}\left[\sum_{i=1}^{N}\left(\left\Vert \boldsymbol{\theta}_{i}\left(1\right)\right\Vert +\left\Vert \boldsymbol{\theta}^{*}\right\Vert \right)^{2}\right]\leq4NR^{2}.
\end{equation}
Thus $d_{t}\leq\Psi\left(\lambda_{\mathsf{I}}+\lambda_{\mathsf{II}}\frac{\alpha_{0}}{\beta_{0}^{2}}t^{-\frac{1}{2}}\right)^{2}t^{-\frac{1}{2}}$
holds for all $t\leq t_{0}$ if 
\begin{equation}
\Psi\geq\max_{t\leq t_{0}}\left(4NR^{2}\left(\lambda_{\mathsf{I}}+\lambda_{\mathsf{II}}\frac{\alpha_{0}}{\beta_{0}^{2}}t^{-\frac{1}{2}}\right)^{-2}t^{\frac{1}{2}}\right)=4NR^{2}\left(\lambda_{\mathsf{I}}+\lambda_{\mathsf{II}}\frac{\alpha_{0}}{\beta_{0}^{2}}t_{0}^{-\frac{1}{2}}\right)^{-2}t_{0}^{\frac{1}{2}}.\label{eq:omega_condi1}
\end{equation}

The next step is to show that, for all $t\geq t_{0}$, $\sqrt{d_{t}}\leq\left(\lambda_{\mathsf{I}}+\lambda_{\mathsf{II}}\frac{\alpha_{0}}{\beta_{0}^{2}}t^{-\frac{1}{2}}\right)t^{-\frac{1}{4}}\sqrt{\Psi}$
leads to $\sqrt{d_{t+1}}\leq\left(\lambda_{\mathsf{I}}+\lambda_{\mathsf{II}}\frac{\alpha_{0}}{\beta_{0}^{2}}\left(t+1\right)^{-\frac{1}{2}}\right)\left(t+1\right)^{-\frac{1}{4}}\sqrt{\Psi}$
if $\Psi$ satisfies (\ref{eq:const_strong}). 

With $\alpha_{t}=\alpha_{0}t^{-1}$, $\beta_{t}=\beta_{0}t^{-\frac{1}{4}}$,
and $\sqrt{d_{t}}\leq\left(\lambda_{\mathsf{I}}+\lambda_{\mathsf{II}}\frac{\alpha_{0}}{\beta_{0}^{2}}t^{-\frac{1}{2}}\right)t^{-\frac{1}{4}}\sqrt{\Psi}$,
we have the following bound according to (\ref{eq:evo_strong}):
\begin{align}
 & d_{t+1}\leq\left(1-2\mu\alpha_{0}t^{-1}\right)d_{t}+\left(\lambda_{\mathsf{I}}+\lambda_{\mathsf{II}}\frac{\alpha_{0}}{\beta_{0}^{2}}t^{-\frac{1}{2}}\right)\alpha_{0}\beta_{0}t^{-\frac{5}{4}}\sqrt{d_{t}}+\lambda_{\mathsf{III}}\frac{\alpha_{0}^{2}}{\beta_{0}^{2}}t^{-\frac{3}{2}}\nonumber \\
 & \leq\!\left(1\!-2\mu\alpha_{0}t^{-1}\right)\!\left(\lambda_{\mathsf{I}}+\lambda_{\mathsf{II}}\frac{\alpha_{0}}{\beta_{0}^{2}}t^{-\frac{1}{2}}\right)^{\!2}\!t^{-\frac{1}{2}}\Psi\!+\!\left(\!\lambda_{\mathsf{I}}\!+\!\lambda_{\mathsf{II}}\frac{\alpha_{0}}{\beta_{0}^{2}}t^{-\frac{1}{2}}\!\right)^{\!2}\!\alpha_{0}\beta_{0}t^{-\frac{3}{2}}\sqrt{\Psi}+\lambda_{\mathsf{III}}\frac{\alpha_{0}^{2}}{\beta_{0}^{2}}t^{-\frac{3}{2}},\label{eq:evo_strong1}
\end{align}
Note that as $t\geq t_{0}=\left\lceil 2\mu\alpha_{0}\right\rceil \geq2\mu\alpha_{0}$,
we have 
\begin{equation}
1-2\mu\alpha_{0}t^{-1}\geq1-2\mu\alpha_{0}\left(2\mu\alpha_{0}\right)^{-1}=0.
\end{equation}
A sufficient condition to have $\sqrt{d_{t+1}}\leq\left(\lambda_{\mathsf{I}}+\lambda_{\mathsf{II}}\frac{\alpha_{0}}{\beta_{0}^{2}}\left(t+1\right)^{-\frac{1}{2}}\right)\left(t+1\right)^{-\frac{1}{4}}\sqrt{\Psi}$
is to have the following satisfied
\begin{align}
 & \left(1-2\mu\alpha_{0}t^{-1}\right)\!\left(\!\lambda_{\mathsf{I}}+\lambda_{\mathsf{II}}\frac{\alpha_{0}}{\beta_{0}^{2}}t^{-\frac{1}{2}}\!\right)^{\!2}\!t^{-\frac{1}{2}}\Psi\!+\!\left(\!\lambda_{\mathsf{I}}+\!\lambda_{\mathsf{II}}\frac{\alpha_{0}}{\beta_{0}^{2}}t^{-\frac{1}{2}}\!\right)^{\!2}\!\!\alpha_{0}\beta_{0}t^{-\frac{3}{2}}\!\sqrt{\Psi}+\!\lambda_{\mathsf{III}}\frac{\alpha_{0}^{2}}{\beta_{0}^{2}}t^{-\frac{3}{2}}\nonumber \\
 & \leq\left(\lambda_{\mathsf{I}}+\lambda_{\mathsf{II}}\frac{\alpha_{0}}{\beta_{0}^{2}}\left(t+1\right)^{-\frac{1}{2}}\right)^{2}\left(t+1\right)^{-\frac{1}{2}}\Psi.\label{eq:evo_strong_2}
\end{align}
In other words, we need to show that there exists $\Psi<\infty$ such
that (\ref{eq:evo_strong_2}) holds for any $t\geq t_{0}$.

We can rewrite (\ref{eq:evo_strong_2}) as, 
\begin{align}
\left(\!t-\!\left(\!\frac{\lambda_{\mathsf{I}}+\lambda_{\mathsf{II}}\frac{\alpha_{0}}{\beta_{0}^{2}}\left(t+1\right)^{-\frac{1}{2}}}{\lambda_{\mathsf{I}}+\lambda_{\mathsf{II}}\frac{\alpha_{0}}{\beta_{0}^{2}}t^{-\frac{1}{2}}}\!\right)^{\!\!2}\!\!\left(t+1\right)^{-\frac{1}{2}}t^{\frac{3}{2}}-2\mu\alpha_{0}\!\right)\!\!\Psi\nonumber \\
+\alpha_{0}\beta_{0}\sqrt{\Psi}+\frac{\lambda_{\mathsf{III}}\frac{\alpha_{0}^{2}}{\beta_{0}^{2}}}{\left(\lambda_{\mathsf{I}}+\lambda_{\mathsf{II}}\frac{\alpha_{0}}{\beta_{0}^{2}}t^{-\frac{1}{2}}\right)^{\!2}}\leq0.\label{eq:omeg_1}
\end{align}
We need to use the bounds
\begin{equation}
\frac{\lambda_{\mathsf{III}}\frac{\alpha_{0}^{2}}{\beta_{0}^{2}}}{\left(\lambda_{\mathsf{I}}+\lambda_{\mathsf{II}}\frac{\alpha_{0}}{\beta_{0}^{2}}t^{-\frac{1}{2}}\right)^{2}}\leq\frac{\lambda_{\mathsf{III}}\alpha_{0}^{2}}{\lambda_{\mathsf{I}}^{2}\beta_{0}^{2}},\label{eq:bound1}
\end{equation}
and
\begin{align}
 & t-\!\left(\!\frac{\lambda_{\mathsf{I}}+\lambda_{\mathsf{II}}\frac{\alpha_{0}}{\beta_{0}^{2}}\left(t+1\right)^{-\frac{1}{2}}}{\lambda_{\mathsf{I}}+\lambda_{\mathsf{II}}\frac{\alpha_{0}}{\beta_{0}^{2}}t^{-\frac{1}{2}}}\!\right)^{\!\!2}\!\!\left(t+1\right)^{-\frac{1}{2}}t^{\frac{3}{2}}\nonumber \\
 & =t-\!\!\left(\!1+\frac{\lambda_{\mathsf{II}}\frac{\alpha_{0}}{\beta_{0}^{2}}\left(\left(t+1\right)^{-\frac{1}{2}}\!-t^{-\frac{1}{2}}\right)}{\lambda_{\mathsf{I}}+\lambda_{\mathsf{II}}\frac{\alpha_{0}}{\beta_{0}^{2}}t^{-\frac{1}{2}}}\!\right)^{\!\!2}\!\!\left(t+1\right)^{-\frac{1}{2}}t^{\frac{3}{2}}\nonumber \\
 & =t-\!\left(\!1+2\lambda_{\mathsf{II}}\frac{\alpha_{0}}{\beta_{0}^{2}}\frac{\left(t+1\right)^{-\frac{1}{2}}-t^{-\frac{1}{2}}}{\lambda_{\mathsf{I}}+\lambda_{\mathsf{II}}\frac{\alpha_{0}}{\beta_{0}^{2}}t^{-\frac{1}{2}}}+\!\left(\!\lambda_{\mathsf{II}}\frac{\alpha_{0}}{\beta_{0}^{2}}\frac{\left(t+1\right)^{-\frac{1}{2}}-t^{-\frac{1}{2}}}{\lambda_{\mathsf{I}}+\lambda_{\mathsf{II}}\frac{\alpha_{0}}{\beta_{0}^{2}}t^{-\frac{1}{2}}}\right)^{\!2}\right)\!\left(t+1\right)^{-\frac{1}{2}}t^{\frac{3}{2}}\nonumber \\
 & \leq\left(t^{-\frac{1}{2}}-\left(t+1\right)^{-\frac{1}{2}}\right)t^{\frac{3}{2}}+2\lambda_{\mathsf{II}}\frac{\alpha_{0}}{\beta_{0}^{2}}\frac{\left(t^{-\frac{1}{2}}-\left(t+1\right)^{-\frac{1}{2}}\right)t^{\frac{3}{2}}}{\lambda_{\mathsf{I}}+\lambda_{\mathsf{II}}\frac{\alpha_{0}}{\beta_{0}^{2}}t^{-\frac{1}{2}}}\left(t+1\right)^{-\frac{1}{2}}\nonumber \\
 & \leq\frac{1}{2}+\lambda_{\mathsf{II}}\frac{\alpha_{0}}{\beta_{0}^{2}}\frac{\left(t+1\right)^{-\frac{1}{2}}}{\lambda_{\mathsf{I}}+\lambda_{\mathsf{II}}\frac{\alpha_{0}}{\beta_{0}^{2}}t^{-\frac{1}{2}}}\label{eq:bound2}\\
 & \leq\frac{1}{2}+\lambda_{\mathsf{II}}\frac{\alpha_{0}}{\beta_{0}^{2}}\frac{\left(t+1\right)^{-\frac{1}{2}}}{\lambda_{\mathsf{II}}\frac{\alpha_{0}}{\beta_{0}^{2}}t^{-\frac{1}{2}}}=\frac{1}{2}+\sqrt{\frac{t}{t+1}}\leq\frac{3}{2},\label{eq:bound3}
\end{align}
where (\ref{eq:bound2}) is by
\begin{align}
 & \left(t^{-\frac{1}{2}}-\left(t+1\right)^{-\frac{1}{2}}\right)t^{\frac{3}{2}}\nonumber \\
 & =\left(t^{-\frac{1}{2}}\left(t+1\right)^{\frac{1}{2}}-1\right)t^{\frac{3}{2}}\left(t+1\right)^{-\frac{1}{2}}=\left(\left(1+\frac{1}{t}\right)^{\frac{1}{2}}-1\right)t^{\frac{3}{2}}\left(t+1\right)^{-\frac{1}{2}}\nonumber \\
 & \leq\frac{1}{2}t^{-1}t^{\frac{3}{2}}\left(t+1\right)^{-\frac{1}{2}}=\frac{1}{2}\left(\frac{t}{t+1}\right)^{\frac{1}{2}}\leq\frac{1}{2}.
\end{align}
Using (\ref{eq:bound1}) and (\ref{eq:bound3}), we get
\begin{align}
 & \left(\!t-\!\left(\!\frac{\lambda_{\mathsf{I}}+\lambda_{\mathsf{II}}\frac{\alpha_{0}}{\beta_{0}^{2}}\left(t+1\right)^{-\frac{1}{2}}}{\lambda_{\mathsf{I}}+\lambda_{\mathsf{II}}\frac{\alpha_{0}}{\beta_{0}^{2}}t^{-\frac{1}{2}}}\!\right)^{\!\!2}\!\!\left(t+1\right)^{-\frac{1}{2}}t^{\frac{3}{2}}-2\mu\alpha_{0}\!\right)\!\Psi+\alpha_{0}\beta_{0}\sqrt{\Psi}+\frac{\lambda_{\mathsf{III}}\frac{\alpha_{0}^{2}}{\beta_{0}^{2}}}{\left(\lambda_{\mathsf{I}}+\lambda_{\mathsf{II}}\frac{\alpha_{0}}{\beta_{0}^{2}}t^{-\frac{1}{2}}\right)^{2}}\nonumber \\
 & \leq\left(\frac{3}{2}-2\mu\alpha_{0}\right)\Psi+\alpha_{0}\beta_{0}\sqrt{\Psi}+\frac{\lambda_{\mathsf{III}}\alpha_{0}^{2}}{\lambda_{\mathsf{I}}^{2}\beta_{0}^{2}}.
\end{align}
Then (\ref{eq:omeg_1}) can be ensured if the following holds
\begin{equation}
\left(-\frac{3}{2}+2\mu\alpha_{0}\right)\Psi-\alpha_{0}\beta_{0}\sqrt{\Psi}-\frac{\lambda_{\mathsf{III}}\alpha_{0}^{2}}{\lambda_{\mathsf{I}}^{2}\beta_{0}^{2}}\geq0.\label{eq:omeg_2}
\end{equation}
As $\sqrt{\Psi}>0$, We can solve (\ref{eq:omeg_2}) to get:
\begin{equation}
\sqrt{\Psi}\geq\frac{\beta_{0}+\sqrt{\beta_{0}^{2}+2\left(4\mu\alpha_{0}-3\right)\frac{\lambda_{\mathsf{III}}}{\lambda_{\mathsf{I}}^{2}\beta_{0}^{2}}}}{4\mu-\frac{3}{\alpha_{0}}}.\label{eq:omega_condi2}
\end{equation}
Therefore, for all $t\geq t_{0}$, $\sqrt{d_{t+1}}\leq\left(\lambda_{\mathsf{I}}+\lambda_{\mathsf{II}}\frac{\alpha_{0}}{\beta_{0}^{2}}\left(t+1\right)^{-\frac{1}{2}}\right)\left(t+1\right)^{-\frac{1}{4}}\sqrt{\Psi}$
can be deduced under condition $\sqrt{d_{t}}\leq\left(\lambda_{\mathsf{I}}+\lambda_{\mathsf{II}}\frac{\alpha_{0}}{\beta_{0}^{2}}t^{-\frac{1}{2}}\right)t^{-\frac{1}{4}}\sqrt{\Psi}$,
as long as (\ref{eq:omega_condi2}) holds.

As a conclusion, we have shown (\ref{eq:rate_strong}) when $\Psi$
satisfies both (\ref{eq:omega_condi1}) and (\ref{eq:omega_condi2}). 

\section{\label{subsec:appx_f}Proof of Lemma~\ref{lem:step_const}}

The proof is by applying the lower bounds of two fundamental functions.
We have
\begin{align}
\frac{\beta_{0}+\sqrt{\beta_{0}^{2}+2\left(4\mu\alpha_{0}-3\right)\frac{\lambda_{\mathsf{III}}}{\lambda_{\mathsf{I}}^{2}\beta_{0}^{2}}}}{4\mu-\frac{3}{\alpha_{0}}} & \geq\frac{3^{\frac{3}{4}}\left(2\left(4\mu\alpha_{0}-3\right)\frac{\lambda_{\mathsf{III}}}{\lambda_{\mathsf{I}}^{2}}\right)^{\frac{1}{4}}}{4\mu-\frac{3}{\alpha_{0}}}\label{eq:omega_3}\\
 & =3^{\frac{3}{4}}\left(2\lambda_{\mathsf{I}}^{-2}\lambda_{\mathsf{III}}\right)^{\frac{1}{4}}\left(4\mu\alpha_{0}-3\right)^{-\frac{3}{4}}\alpha_{0}\nonumber \\
 & \geq3^{\frac{3}{4}}\left(2\lambda_{\mathsf{I}}^{-2}\lambda_{\mathsf{III}}\right)^{\frac{1}{4}}3^{-\frac{1}{2}}\mu^{-1}=\left(6\lambda_{\mathsf{I}}^{-2}\lambda_{\mathsf{III}}\right)^{\frac{1}{4}}\mu^{-1}\label{eq:omega_5}
\end{align}
where (\ref{eq:omega_3}) is obtained by fact that for any $a>0$
and $\beta_{0}>0$, the function $h_{1}\left(\beta_{0}\right)=\beta_{0}+\sqrt{\beta_{0}^{2}+\vartheta\beta_{0}^{-2}}$
takes its minimum value when $\beta_{0}^{*}=3^{-\frac{1}{4}}\vartheta^{\frac{1}{4}}$
and 
\[
h_{1}\left(\beta_{0}\right)\geq h_{1}\left(\beta_{0}^{*}\right)=\left(3^{-\frac{1}{4}}+\sqrt{3^{\frac{1}{2}}+3^{-\frac{1}{2}}}\right)\vartheta^{\frac{1}{4}}=3^{\frac{3}{4}}\vartheta^{\frac{1}{4}},
\]
here $\vartheta=2\left(4\mu\alpha_{0}-3\right)\frac{\lambda_{\mathsf{III}}}{\lambda_{\mathsf{I}}^{2}}$
and the equality of (\ref{eq:omega_3}) holds when 
\begin{equation}
\beta_{0}^{*}=\left(\frac{2\left(4\mu\alpha_{0}-3\right)\lambda_{\mathsf{III}}}{3\lambda_{\mathsf{I}}^{2}}\right)^{\frac{1}{4}};\label{eq:omega_4}
\end{equation}
(\ref{eq:omega_5}) comes from the property of another function $h_{2}\left(\alpha_{0}\right)=\left(4\mu\alpha_{0}-3\right)^{-\frac{3}{4}}\alpha_{0}$
with $\alpha_{0}>\frac{3}{4\mu}$, it is straightforward to show that
as 
\begin{equation}
\alpha_{0}^{*}=\frac{3}{\mu},\label{eq:omega_6}
\end{equation}
$h_{2}$ takes its minimum value $h_{2}\left(\alpha_{0}^{*}\right)=\left(4\mu\frac{3}{\mu}-3\right)^{-\frac{3}{4}}\frac{3}{\mu}=3^{-\frac{1}{2}}\mu^{-1}$.
Substituting (\ref{eq:omega_6}) into (\ref{eq:omega_4}), we get
\[
\beta_{0}^{*}=\left(6\lambda_{\mathsf{I}}^{-2}\lambda_{\mathsf{III}}\right)^{\frac{1}{4}}.
\]

\section{\label{subsec:appx_g}Proof of Lemma~\ref{lem:general_lem1}}

We start with the definition of $\overline{\boldsymbol{\theta}}^{\left(T\right)}$.
Similar to (\ref{eq:p1}), we have

\begin{align}
 & f\!\left(\overline{\boldsymbol{\theta}}^{\left(T\right)}\right)=f\!\left(\!\frac{1}{NT}\sum_{t=1}^{T}\sum_{i=1}^{N}\boldsymbol{\theta}_{i}\left(t\right)\!\right)\!\leq\frac{1}{T}\sum_{t=1}^{T}f\!\left(\!\frac{1}{N}\sum_{i=1}^{N}\boldsymbol{\theta}_{i}\left(t\right)\!\right)=\frac{1}{T}\sum_{t=1}^{T}f\left(\overline{\boldsymbol{\theta}}\left(t\right)\right)\nonumber \\
 & =\frac{1}{NT}\sum_{t=1}^{T}\sum_{i=1}^{N}f_{i}\left(\overline{\boldsymbol{\theta}}\left(t\right)\right)\leq\frac{1}{NT}\sum_{t=1}^{T}\sum_{i=1}^{N}\left(f_{i}\left(\boldsymbol{\theta}_{i}\left(t\right)\right)+\left|f_{i}\left(\overline{\boldsymbol{\theta}}\left(t\right)\right)-f_{i}\left(\boldsymbol{\theta}_{i}\left(t\right)\right)\right|\right)\nonumber \\
 & \leq\frac{1}{NT}\sum_{t=1}^{T}\sum_{i=1}^{N}\left(f_{i}\left(\boldsymbol{\theta}_{i}\left(t\right)\right)+\ell_{i}\left\Vert \overline{\boldsymbol{\theta}}\left(t\right)-\boldsymbol{\theta}_{i}\left(t\right)\right\Vert \right)\leq\frac{1}{NT}\sum_{t=1}^{T}\sum_{i=1}^{N}f_{i}\left(\boldsymbol{\theta}_{i}\left(t\right)\right)+\ell\delta\label{eq:general__}
\end{align}
which is by the assumption that $f_{i}$ is $\ell_{i}$-Lipschitz
(A5). Recall that $\ell=\frac{1}{N}\sum_{i=1}^{N}\ell_{i}$. We have
also used the bound $\left\Vert \overline{\boldsymbol{\theta}}\left(t\right)-\boldsymbol{\theta}_{i}\left(t\right)\right\Vert \leq\delta$
that is presented in Lemma~\ref{lem:dist_aver}. 

From (\ref{eq:general__}), we can deduce that
\begin{align}
 & \mathbb{E}\left[f\left(\overline{\boldsymbol{\theta}}^{\left(T\right)}\right)-f\left(\boldsymbol{\theta}^{*}\right)\right]\leq\frac{1}{NT}\sum_{t=1}^{T}\sum_{i=1}^{N}\left(\mathbb{E}\left[f_{i}\left(\boldsymbol{\theta}_{i}\left(t\right)\right)\right]-f_{i}\left(\boldsymbol{\theta}^{*}\right)\right)+\frac{1}{N}\sum_{i=1}^{N}\ell_{i}\delta\nonumber \\
 & =\frac{1}{NT}\sum_{t=1}^{T}\sum_{i=1}^{N}\left(\mathbb{E}\left[\widetilde{f}_{i}\left(\boldsymbol{\theta}_{i}\left(t\right)\right)\right]-\widetilde{f}_{i}\left(\widetilde{\boldsymbol{\theta}}^{*}\right)\right)+\frac{1}{NT}\sum_{t=1}^{T}\sum_{i=1}^{N}\left(\widetilde{f}_{i}\left(\widetilde{\boldsymbol{\theta}}^{*}\right)-f_{i}\left(\boldsymbol{\theta}^{*}\right)\right)\nonumber \\
 & \qquad+\frac{1}{NT}\sum_{t=1}^{T}\sum_{i=1}^{N}\mathbb{E}\left[f_{i}\left(\boldsymbol{\theta}_{i}\left(t\right)\right)-\widetilde{f}_{i}\left(\boldsymbol{\theta}_{i}\left(t\right)\right)\right]+\ell\delta.\label{eq:genral_error}
\end{align}
Clearly, our next step is to find upper bounds of $f_{i}\left(\boldsymbol{\theta}_{i}\left(t\right)\right)-\widetilde{f}_{i}\left(\boldsymbol{\theta}_{i}\left(t\right)\right)$and
of $\widetilde{f}_{i}\left(\widetilde{\boldsymbol{\theta}}^{*}\right)-f_{i}\left(\boldsymbol{\theta}^{*}\right)$.
By definition of $\widetilde{f}_{i}$, we know that 
\begin{align}
\widetilde{f}_{i}\left(\boldsymbol{\theta}_{i}\left(t\right)\right) & =\mathbb{E}_{\boldsymbol{\varpi}\in\mathbb{R}^{M}:\left\Vert \boldsymbol{\varpi}\right\Vert \leq1}\left[f_{i}\left(\boldsymbol{\theta}_{i}\left(t\right)+\beta\boldsymbol{\varpi}\right)\right]\nonumber \\
 & \geq\mathbb{E}_{\boldsymbol{\varpi}\in\mathbb{R}^{M}:\left\Vert \boldsymbol{\varpi}\right\Vert \leq1}\left[f_{i}\left(\boldsymbol{\theta}_{i}\left(t\right)\right)-\left|f_{i}\left(\boldsymbol{\theta}_{i}\left(t\right)+\beta\boldsymbol{\varpi}\right)-f_{i}\left(\boldsymbol{\theta}_{i}\left(t\right)\right)\right|\right]\nonumber \\
 & \geq f_{i}\left(\boldsymbol{\theta}_{i}\left(t\right)\right)-\mathbb{E}_{\boldsymbol{\varpi}\in\mathbb{R}^{M}:\left\Vert \boldsymbol{\varpi}\right\Vert \leq1}\left[\ell_{i}\left\Vert \beta\boldsymbol{\varpi}\right\Vert \right]\geq f_{i}\left(\boldsymbol{\theta}_{i}\left(t\right)\right)-\beta\ell_{i},\label{eq:gap_smooth_1}
\end{align}
leading to
\begin{equation}
\frac{1}{NT}\sum_{t=1}^{T}\sum_{i=1}^{N}\mathbb{E}\left[f_{i}\left(\boldsymbol{\theta}_{i}\left(t\right)\right)-\widetilde{f}_{i}\left(\boldsymbol{\theta}_{i}\left(t\right)\right)\right]\leq\ell\beta.\label{eq:gap_smooth}
\end{equation}
Similar to (\ref{eq:gap_smooth_1}), we can also deduce that $\widetilde{f}_{i}\left(\boldsymbol{\theta}\right)\leq f_{i}\left(\boldsymbol{\theta}\right)+\ell_{i}\beta$,
which implies 
\[
\frac{1}{N}\sum_{i\in\mathcal{N}}\widetilde{f}_{i}\left(\widetilde{\boldsymbol{\theta}}^{*}\right)=\min_{\boldsymbol{\theta}\in\widetilde{\mathcal{K}}}\frac{1}{N}\sum_{i\in\mathcal{N}}\widetilde{f}_{i}\left(\boldsymbol{\theta}\right)\leq\min_{\boldsymbol{\theta}\in\widetilde{\mathcal{K}}}\frac{1}{N}\sum_{i\in\mathcal{N}}f_{i}\left(\boldsymbol{\theta}\right)+\ell\beta.
\]
By assuming that $\boldsymbol{\theta}^{*}\in\widetilde{\mathcal{K}}$,
we also have $f\left(\boldsymbol{\theta}^{*}\right)=\min_{\boldsymbol{\theta}\in\widetilde{\mathcal{K}}}\frac{1}{N}\sum_{i\in\mathcal{N}}f_{i}\left(\boldsymbol{\theta}\right)$.
Thus 
\begin{equation}
\frac{1}{NT}\sum_{t=1}^{T}\sum_{i=1}^{N}\left(\widetilde{f}_{i}\left(\widetilde{\boldsymbol{\theta}}^{*}\right)-f_{i}\left(\boldsymbol{\theta}^{*}\right)\right)\leq\ell\beta.\label{eq:gap_opt}
\end{equation}
We can introduce (\ref{eq:gap_smooth}) and (\ref{eq:gap_opt}) into
(\ref{eq:genral_error}) to obtain (\ref{eq:genral_error-1}), which
concludes the proof. 

\section{\label{subsec:appx_h}Proof of Lemma~\ref{lem:general_lem2}}

Denote 
\[
\widetilde{D}_{t}=\sum_{i=1}^{N}\left\Vert \boldsymbol{\theta}_{i}\left(t\right)-\widetilde{\boldsymbol{\theta}}^{*}\right\Vert ^{2}
\]
 Using the similar steps as in (\ref{eq:dt1-dt}), (\ref{eq:xg3}),
and (\ref{eq:xg4}), we can get

\begin{align}
 & \widetilde{D}_{t+1}=\sum_{i=1}^{N}\left\Vert \boldsymbol{\theta}_{i}\left(t+1\right)-\widetilde{\boldsymbol{\theta}}^{*}\right\Vert ^{2}\nonumber \\
 & =\sum_{i=1}^{N}\left\Vert \overline{\boldsymbol{\theta}}\left(t\right)-\widetilde{\boldsymbol{\theta}}^{*}\right\Vert ^{2}+\sum_{i=1}^{N}\left\Vert \sum_{j=1}^{N}A_{i,j}\left(t\right)\boldsymbol{\theta}_{j}\left(t\right)-\overline{\boldsymbol{\theta}}\left(t\right)-\alpha\widehat{\boldsymbol{g}}_{i}\left(t\right)\right\Vert ^{2}\nonumber \\
 & -2\alpha\sum_{i=1}^{N}\left\langle \overline{\boldsymbol{\theta}}\left(t\right)-\widetilde{\boldsymbol{\theta}}^{*},\widehat{\boldsymbol{g}}_{i}\left(t\right)\right\rangle +2\sum_{i=1}^{N}\left\langle \overline{\boldsymbol{\theta}}\left(t\right)-\widetilde{\boldsymbol{\theta}}^{*},\sum_{j=1}^{N}A_{i,j}\left(t\right)\boldsymbol{\theta}_{j}\left(t\right)-\overline{\boldsymbol{\theta}}\left(t\right)\right\rangle \nonumber \\
 & \leq\widetilde{D}_{t}+2N\delta^{2}+2NC^{2}\frac{\alpha^{2}}{\beta^{2}}-2\alpha\sum_{i=1}^{N}\left\langle \overline{\boldsymbol{\theta}}\left(t\right)-\widetilde{\boldsymbol{\theta}}^{*},\widehat{\boldsymbol{g}}_{i}\left(t\right)\right\rangle .\label{eq:Dapprox}
\end{align}
Note that in this part of analysis, the step-sizes are constant $\alpha$
and $\beta$, besides the perturbation vector is a unit vector so
that $\left\Vert \widehat{\boldsymbol{g}}_{i}\left(t\right)\right\Vert \leq C/\beta$. 

Now we need to find a lower bound of $\sum_{i=1}^{N}\left\langle \overline{\boldsymbol{\theta}}\left(t\right)-\widetilde{\boldsymbol{\theta}}^{*},\widehat{\boldsymbol{g}}_{i}\left(t\right)\right\rangle $.
Again, we introduce the stochastic noise $\boldsymbol{e}_{i}\left(t\right)$
such that $\boldsymbol{e}_{i}\left(t\right)=\widehat{\boldsymbol{g}}_{i}\left(t\right)-\mathbb{E}\left[\widehat{\boldsymbol{g}}_{i}\left(t\right)\mid\boldsymbol{\theta}_{i}\left(t\right)\right]$.
We have
\begin{align}
 & \sum_{i=1}^{N}\left\langle \overline{\boldsymbol{\theta}}\left(t\right)-\widetilde{\boldsymbol{\theta}}^{*},\widehat{\boldsymbol{g}}_{i}\left(t\right)-\boldsymbol{e}_{i}\left(t\right)\right\rangle =\sum_{i=1}^{N}\left\langle \overline{\boldsymbol{\theta}}\left(t\right)-\widetilde{\boldsymbol{\theta}}^{*},\mathbb{E}\left[\widehat{\boldsymbol{g}}_{i}\left(t\right)\mid\boldsymbol{\theta}_{i}\left(t\right)\right]\right\rangle \nonumber \\
 & =\sum_{i=1}^{N}\left\langle \overline{\boldsymbol{\theta}}\left(t\right)-\boldsymbol{\theta}_{i}\left(t\right)+\boldsymbol{\theta}_{i}-\widetilde{\boldsymbol{\theta}}^{*},\mathbb{E}\left[\widehat{\boldsymbol{g}}_{i}\left(t\right)\mid\boldsymbol{\theta}_{i}\left(t\right)\right]\right\rangle \nonumber \\
 & =\sum_{i=1}^{N}\!\left\langle \!\!\boldsymbol{\theta}_{i}\!\left(t\right)\!-\widetilde{\boldsymbol{\theta}}^{*},\frac{\nabla\widetilde{f}_{i}\left(\boldsymbol{\theta}_{i}\!\left(t\right)\right)}{M}\!\right\rangle \!+\sum_{i=1}^{N}\!\left\langle \overline{\boldsymbol{\theta}}\!\left(t\right)-\boldsymbol{\theta}_{i}\!\left(t\right),\mathbb{E}\left[\widehat{\boldsymbol{g}}_{i}\!\left(t\right)\mid\boldsymbol{\theta}_{i}\!\left(t\right)\right]\right\rangle \label{eq:g1}\\
 & \geq\frac{1}{M}\sum_{i=1}^{N}\left(\widetilde{f}_{i}\left(\boldsymbol{\theta}_{i}\!\left(t\right)\right)-\widetilde{f}_{i}\left(\widetilde{\boldsymbol{\theta}}^{*}\right)\right)-\sum_{i=1}^{N}\left\Vert \overline{\boldsymbol{\theta}}\!\left(t\right)-\boldsymbol{\theta}_{i}\!\left(t\right)\right\Vert \left\Vert \mathbb{E}\left[\widehat{\boldsymbol{g}}_{i}\!\left(t\right)\mid\boldsymbol{\theta}_{i}\!\left(t\right)\right]\right\Vert \label{eq:g2}\\
 & \geq\frac{1}{M}\sum_{i=1}^{N}\left(\widetilde{f}_{i}\left(\boldsymbol{\theta}_{i}\left(t\right)\right)-\widetilde{f}_{i}\left(\widetilde{\boldsymbol{\theta}}^{*}\right)\right)-NC\frac{\delta}{\beta},\label{eq:g3}
\end{align}
in which (\ref{eq:g1}) can be proved by applying Lemma~\ref{lem:bias_general},
i.e., 
\begin{align*}
\mathbb{E}\left[\widehat{\boldsymbol{g}}_{i}\left(t\right)\!\mid\!\boldsymbol{\theta}_{i}\!\left(t\right)\right] & =\mathbb{E}_{\boldsymbol{\nu}_{i,t},\boldsymbol{\xi}_{i,t}}\!\!\left[\frac{\boldsymbol{\nu}_{i,t}F_{i}\left(\boldsymbol{\theta}_{i}\!\left(t\right)+\beta\boldsymbol{\nu}_{i,t};\boldsymbol{\xi}_{i,t}\right)}{\beta}\right]\!=\frac{1}{\beta}\mathbb{E}_{\boldsymbol{\nu}_{i,t}}\left[\boldsymbol{\nu}_{i,t}f_{i}\!\left(\boldsymbol{\theta}_{i}\!\left(t\right)+\beta_{t}\boldsymbol{\nu}_{i,t}\right)\right]\\
 & =\frac{1}{\beta}\frac{\beta}{M}\nabla\widetilde{f}_{i}\left(\boldsymbol{\theta}_{i}\left(t\right)\right)=\frac{1}{M}\nabla\widetilde{f}_{i}\left(\boldsymbol{\theta}_{i}\left(t\right)\right);
\end{align*}
(\ref{eq:g2}) comes from the fact that $\widetilde{f}_{i}\left(\boldsymbol{\theta}_{i}\left(t\right)\right)$
is a convex function of $\boldsymbol{\theta}_{i}\left(t\right)$;
we use Lemma~\ref{lem:dist_aver} to get $\left\Vert \overline{\boldsymbol{\theta}}\left(t\right)-\boldsymbol{\theta}_{i}\left(t\right)\right\Vert \leq\delta$
and also the fact that $\left\Vert \widehat{\boldsymbol{g}}_{i}\left(t\right)\right\Vert \leq C/\beta$,
thus (\ref{eq:g3}) can be obtained.

By substituting (\ref{lem:dist_aver}) into (\ref{eq:Dapprox}), we
have 
\begin{align}
\widetilde{D}_{t+1} & \leq\widetilde{D}_{t}+2N\delta^{2}+2NC^{2}\frac{\alpha^{2}}{\beta^{2}}-2\alpha\frac{1}{M}\sum_{i=1}^{N}\left(\widetilde{f}_{i}\left(\boldsymbol{\theta}_{i}\left(t\right)\right)-\widetilde{f}_{i}\left(\widetilde{\boldsymbol{\theta}}^{*}\right)\right)\nonumber \\
 & \qquad+2\alpha NC\frac{\delta}{\beta}-2\alpha\sum_{i=1}^{N}\left\langle \overline{\boldsymbol{\theta}}\left(t\right)-\widetilde{\boldsymbol{\theta}}^{*},\boldsymbol{e}_{i}\left(t\right)\right\rangle .\label{eq:Dapprox_2}
\end{align}
Denote $\widetilde{d}_{t}=\mathbb{E}[\widetilde{D}_{t}]$, we can
take expectation on both sides of (\ref{eq:Dapprox_2}) to get
\begin{align}
\sum_{i=1}^{N}\left(\mathbb{E}\left[\widetilde{f}_{i}\left(\boldsymbol{\theta}_{i}\left(t\right)\right)\right]-\widetilde{f}_{i}\left(\widetilde{\boldsymbol{\theta}}^{*}\right)\right) & \leq\frac{M}{2\alpha}\left(\widetilde{d}_{t}-\widetilde{d}_{t+1}\right)+NM\frac{\delta^{2}}{\alpha}+NMC^{2}\frac{\alpha}{\beta^{2}}+NMC\frac{\delta}{\beta}.\label{eq:Dapprox3}
\end{align}
We perform the summation of (\ref{eq:Dapprox3}) from $t=1$ to $t=T$:
\begin{align*}
 & \frac{1}{NT}\sum_{t=1}^{T}\sum_{i=1}^{N}\left(\mathbb{E}\left[\widetilde{f}_{i}\left(\boldsymbol{\theta}_{i}\left(t\right)\right)\right]-\widetilde{f}_{i}\left(\widetilde{\boldsymbol{\theta}}^{*}\right)\right)\\
 & \leq\frac{M}{2\alpha}\frac{1}{T}\sum_{t=1}^{T}\left(\widetilde{d}_{t}-\widetilde{d}_{t+1}\right)+NM\frac{\delta^{2}}{\alpha}+NMC^{2}\frac{\alpha}{\beta^{2}}+NMC\frac{\delta}{\beta}\\
 & \leq\frac{M}{2\alpha}\frac{1}{T}\widetilde{d}_{1}+NM\frac{\delta^{2}}{\alpha}+NMC^{2}\frac{\alpha}{\beta^{2}}+NMC\frac{\delta}{\beta}\leq\frac{2MR^{2}}{\alpha T}+M\frac{\delta^{2}}{\alpha}+MC^{2}\frac{\alpha}{\beta^{2}}+MC\frac{\delta}{\beta},
\end{align*}
which concludes the proof.
\bibliographystyle{siamplain}
\bibliography{BiblioWenjie}

\end{document}